\documentclass[onefignum,onetabnum]{article}

\usepackage{mathtools}
\usepackage{amsfonts}
\usepackage{amssymb}
\usepackage{bm}
\usepackage{amsthm}
\usepackage{bbm}
\usepackage{subcaption}
\usepackage{hyperref}
\usepackage{enumitem}
\usepackage{algorithm}
\usepackage{algpseudocode}
\usepackage{mycommands}
\usepackage{customcolors}
\usepackage[]{graphicx}
\usepackage[parfill]{parskip}
\usepackage{tikz}
\usepackage{comment}
\usepackage{pgfplotstable}
\usepackage{pgfplots}
\usepackage{ifthen}
\usepackage{hhline}
\usepackage{capt-of}
\usepackage{xr}

\title{Preconditioning across parameter space for the parametric Helmholtz equation}

\author{Wouter Gerrit van Harten\thanks{IMAPP, Radboud University, The Netherlands
    (\email{w.vanharten@science.ru.nl})}
\and Laura Scarabosio\footnotemark[2]}

\newcommand{\y}{{\bm{y}}}
\newcommand{\ys}{{\hat{\bm{y}}_\ast}}
\newcommand{\yb}{\bar{\bm{y}}}
\newcommand{\yh}{\hat{\bm{y}}}

\newcommand{\yi}{{\bm{y}_i}}

\newcommand{\eq}{=}

\renewcommand{\u}{{\bm{u}}}
\newcommand{\x}{{\bm{x}}}
\newcommand{\dx}{\mathrm{d}\x}
\newcommand{\D}{\mathop{}\!\mathrm{D}_x}

\DeclareMathOperator*{\argmin}{arg\,min}
\DeclareMathOperator*{\argmax}{arg\,max}

\DeclareMathOperator*{\diam}{diam}

\newtheorem{thm}{Theorem}
\newtheorem{assumption}[thm]{Assumption}

\newtheorem{remark}[thm]{Remark}

\algnewcommand{\algorithmicand}{\textbf{ and }}
\algnewcommand{\algorithmicor}{\textbf{ or }}
\algnewcommand{\OR}{\algorithmicor}
\algnewcommand{\AND}{\algorithmicand}
\algnewcommand{\var}{\texttt}

\newcommand*{\email}[1]{\href{mailto:#1}{\nolinkurl{#1}} }

\newenvironment{revenv}{\color{black}}{\ignorespacesafterend\color{black}}
\newenvironment{newenv}{\color{black}}{\ignorespacesafterend\color{black}}

\begin{document}
    \maketitle
    \begin{abstract}
        In this work, we address the efficient computation of parameterized systems of linear equations, with possible nonlinear parameter dependence.
        When the matrix is highly sensitive to the parameters, mean-based preconditioning might not be enough.
        For this scenario, we explore an approach in which several preconditioners are placed in the parameter space during a precomputation step.
        To determine the optimal placement of a limited number of preconditioners, we estimate the expected number of iterations with respect to a given preconditioner a priori and use a location-allocation strategy to optimize the placement of the preconditioners.
        We elaborate on our methodology for the Helmholtz problem with exterior Dirichlet scattering at high frequencies, and we estimate the expected number of GMRES iterations via a gray-box Gaussian process regression approach.
        We illustrate our approach in two practical applications: scattering in a domain with a parametric refractive index and scattering from a scatterer with parameterized shape.
        Using these numerical examples, we show how our methods leads to runtime savings of about an order of magnitude.
        Moreover, we investigate the effect of the parameter dimension and the importance of dimension anisotropy on their efficacy.
    \end{abstract}




    \section{Introduction}\label{sec:introduction}
    In many scientific and engineering applications, solving large systems of parameterized linear equations is a fundamental task.
These systems appear in, among others, uncertainty quantification, PDE-constrained optimization, and inverse problems, where the repeated solution of linear systems is a major computational bottleneck.

For a given linear system, i.e.\ a fixed parameter value, various methods for its efficient solution exist.
While direct methods are feasible for systems of moderate size, larger systems require iterative solvers like Krylov methods~\cite{ipsen1998}.
Their performance is highly dependent on the system's condition number, which can be reduced by using a preconditioner~\cite{wathen2015}.
Performing preconditioning requires itself some computational resources, which are compensated by the realized reduction in Krylov iterations.
Choosing the right preconditioner is application-specific and extensively studied~\cite{pearson2020,wathen2015}.

When solving multiple linear systems, we would like to reduce the computational load.
One could, for example, try to parameterize the preconditioner~\cite{contreras2018}.
\begin{revenv}
    In a somehow similar spirit, a consistent body of recent research has leveraged machine learning tools, mainly neural networks, to improve the performance and generalization capabilities of classical solvers for parameter-dependent systems, e.g.,~\cite{chen2022,luz2020}.
    For the Helmholtz application we will address, this has been considered in~\cite{azulay2022multigrid, cui2025, giraud:hal-05157038, lerer2024multigrid}, and~\cite{khoo2019switchnet} constructs a neural network solver inspired by a preconditioner.
\end{revenv}

Another well-studied option is to recycle either the Krylov subspaces~\cite{parks2006} or the preconditioner(s).
For the Helmholtz application we will consider, Krylov subspace recycling has been considered in~\cite{jin2009}.
In this work, we will focus our attention on preconditioner recycling.
\begin{revenv}
    The strategy we will devise can thus be considered to be complementary to those in~\cite{jin2009,parks2006} and~\cite{azulay2022multigrid,chen2022,contreras2018,cui2025,giraud:hal-05157038,lerer2024multigrid,luz2020}
\end{revenv}

For this, we could take inspiration from~\cite{graham2021}, and apply a single preconditioner to multiple systems.
Initial exploration of this approach has been performed for stochastic spectral finite element methods applied to the diffusion equation~\cite{ghanem1996,keese2004,pellissetti2000} and the Helmholtz equation~\cite{jin2009, wang2019}.
Preconditioning was performed through block-diagonal preconditioning, where a single mean-based preconditioner is used in a block-diagonal setting.
These ideas were then applied to the much more challenging Navier-Stokes equations, where the stochastic Galerkin discretization of said equations was investigated numerically~\cite{powell2012}.
Outside the stochastic Galerkin framework, we can make use of mean-based preconditioning as well, which was initially explored for the diffusion equation~\cite{eiermann2007,ernst2009}.
To get a complete overview of the history of mean-based preconditioning, we refer the reader to an excellent overview by Owen Pembery~\cite[Section~4.7]{pembery2020}.
While mean-based preconditioning is computationally efficient, it might break down if the effect of parameter changes on the matrix is large.

Improvements could be gained by considering multiple preconditioners throughout the parameter space.
This has been investigated in~\cite{venkovic2023} for the diffusion equation with short correlation lengths, focused on Voronoi quantizers placed without prior knowledge of the parameter locations.
For the Helmholtz equation, a greedy approach for known parameter locations was presented in~\cite[Section~4.6]{pembery2020}.
However, this approach requires prior knowledge about the maximum iterations and the distance function on the parameter space, and it might struggle if the problem has hidden parameter anisotropy or a high dimension.

Our approach exploits the structure and a-priori knowledge of parameter locations and overcomes the limitations in~\cite{pembery2020}.
The main novelty of our work is to train a surrogate model to predict the number of Krylov iterations needed and
\begin{revenv}
    parametrically place preconditioners in the parameter space
\end{revenv} while incorporating dimension anisotropy.
Since computing the best locations to place the preconditioners is known to be NP-hard~\cite{sherali1988}, we use a location-allocation heuristic~\cite{brimberg2008}.
After having introduced the problem in Section~\ref{sec:problem-statement}, Section~\ref{sec:placement} explains this process and Section~\ref{sec:estimating-m} discusses methods to obtain such a surrogate using known convergence bounds of the Krylov method.

    \section{Problem statement} \label{sec:problem-statement}
    Let $\{\mathbb{A}(\y_i)\}_{\y_i \in W}$ be a family of parameterized matrices with parameter vectors $W=\{\y_i \in Y\}_{i\in I}$ for some finite index set $I$ and $Y=\prod_{i=1}^N Y_i$ with $Y_i$ a bounded interval in $\mathbb{R}$.
Our goal is to solve the system $\mathbb{A}(\y_i)\u(\y_i)=b$ for all parameter vectors $\y_i$.
We turn our attention to preconditioned Krylov methods and consider the preconditioned system to find the high-dimensional vector $\u(\yi)$ for $\yi\in W$:
\begin{equation}
    \mathbb{P}(\yi)\mathbb{A}(\yi)\u(\yi)=\mathbb{P}(\yi)b.\label{eq:preconditionedsystem}
\end{equation}
We do not specify the exact choice of preconditioner for now, but we will be more explicit in Section~\ref{sec:estimating-m}.
We place $N_{pc}$ preconditioners $\mathbb{P}\left(\hat{\y}_k\right) \in Y$.
These preconditioners form a partition $W_k$ of $W$, each with a selected preconditioner location $\yh_k$.
Then, we set $\mathbb{P}(\yi) = \mathbb{P}(\yh_k)$ for $\yi\in W_k$.

We note that the parameter space $Y$ may be isotropic or anisotropic, implying that the effect of different parameter dimensions on the conditioning of the system may not be equal.
Moreover, we will consider left preconditioning.
Right- and split methods can be treated similarly.

To find a good preconditioning strategy, we will place several preconditioners in the parameter space in the next section.

    \section{Preconditioner placement} \label{sec:placement}
    To place the preconditioners in $Y$, we rely on a problem-specific approximation of the expected number of iterations, which we denote by $\tilde{m}$:
\begin{equation*}
    \tilde{m}(\y,\hat{\y}) \coloneqq \# \text{ Krylov iterations to solve }\mathbb{P}(\yh)\mathbb{A}(\y)\bm{u}=\mathbb{P}(\yh)b, \text{ for }\y,\yh \in Y,
\end{equation*}
and by
\begin{equation}
    N_{ratio} = \frac{\tau_{pc}}{\tau_{Krylov}}\label{eq:Nratio}
\end{equation}
the ratio between the computation time of a preconditioner $\tau_{pc}$ and $\tau_{Krylov}$, the time to perform a single Krylov iteration.
Moreover, we assume that the number of Krylov iterations required is invariant under translations:
\begin{assumption}\label{ass:shift_invariant}
    We have $\tilde{m}(\y,\hat{\y}) = \tilde{m}(\y+\tilde{\y},\hat{\y}+\tilde{\y})$, for all $\y, \tilde{\y},$ and $\hat{\y}$ in $Y$.
\end{assumption}
This way, it is sufficient to define $m(\cdot)$ such that, for $\y,\yh\in Y$,
\begin{equation*}
    m(\y - \hat{\y}) \coloneqq  \tilde{m}(\y,\hat{\y}).
\end{equation*}
\begin{revenv}
    We emphasize that Assumption~\ref{ass:shift_invariant} is a \emph{modelling assumption} which, in general, is not satisfied exactly.
    However, it provides an effective framework for building a good surrogate of $m(\cdot)$ that is efficient to evaluate.
    We will further comment on this assumption in Section~SM1 of the supplementary material.
\end{revenv}
When not fulfilled, we can think of our surrogate as incorporating a modeling approximation.

We must determine $N_{pc}$, the number of preconditioners first.
A small $N_{pc}$ wastes computational potential, while a larger $N_{pc}$ increases the cost of the preconditioning strategy.
For now, we fix $N_{pc}$ and discuss its choice in Section~\ref{subsec:determining-Npc}.
For the time being, we assume we have $m(\cdot)$ and discuss its construction in Section~\ref{sec:estimating-m}.

To place the preconditioners efficiently, we need to solve
\begin{equation}
    \min_{\{\hat{\y}_k\in Y\}_{k=1}^{N_{pc}}} \sum_{\y_j\in W} m\left(\y_j-\hat{\y}\left(\y_j\big|\{\hat{\y}_k\}_{k=1}^{N_{pc}}\right)\right),\label{eq:simplified_optim_problem}
\end{equation}
where 
\begin{equation*}
    \hat{\y}\left(\y_j\big|\{\hat{\y}_k\}_{k=1}^{N_{pc}}\right)=\argmin_{\yh\in \{\hat{\y}_k\}_{k=1}^{N_{pc}}} m(\yi-\yh)
\end{equation*}
maps each parameter location $\yi$ to its assigned preconditioner location $\yh_k$.
Equation~\eqref{eq:simplified_optim_problem} is equivalent to the NP-hard \emph{uncapacitated facility location problem}~\cite{brimberg2008,sherali1988}.
Therefore, we resort to iterative algorithms for computational feasibility, despite the possible occurrence of local minima~\cite{church2022}.

We turn our attention to the well-known \emph{location-allocation} algorithm~\cite{cooper1963,lara2018}, which iterates through a location and an allocation step.
First, the preconditioners are initialized (Section~\ref{subsec:initialization}).
Then the values in $W$ are clustered iteratively, where each step assigns parameters to the preconditioner with the lowest expected number of Krylov iterations.
The preconditioner location is then optimized to minimize Krylov iterations given the assigned parameter values, and this process continues until convergence or diminishing returns (Section~\ref{subsec:location}).

\begin{revenv}
    \begin{remark}\label{rem:kmeans}
        Next to the location-allocation approach discussed here, one could use, for example, $K$-means or competitive learning as suggested in~\cite{venkovic2023}.
        However, these algorithms minimize the square distance to the preconditioner, instead of the expected number of iterations. 
        Therefore, these methods will likely result in a more expensive strategy.
        We investigate this further in Section~SM4 of the supplementary material.
    \end{remark}
\end{revenv}

\subsection{Initialization}\label{subsec:initialization}
Since the optimization problem~\eqref{eq:simplified_optim_problem} has many local minima, good initialization is key.
Several initialization methods have been considered; random initialization proved inefficient, Quasi Monte Carlo sequences can incorporate anisotropy~\cite{howell2009} but require many optimization steps, and Bayesian optimization is effective but suffers from the curse of dimensionality.
We use greedy initialization instead, which is effective and computationally inexpensive.
This means iteratively placing preconditioners at the parameter location with the highest expected number of Krylov iterations.

\subsection{Location-allocation}\label{subsec:location}
After initialization, we iterate through location and allocation steps.
In each location step, we need to find the center of each cell in the partition $W_k$.
This boils down to computing, for every element of the partition:
\begin{equation}
    \argmin_{\ys\in Y} \sum_{\y_j \in W_k} m(\y_j -  \ys). \label{eq:location_step_problem}
\end{equation}
This minimization problem, also known as the \emph{Weber problem} or the \emph{generalized Fermat problem}, is very challenging to solve~\cite{brimberg2008,kalczynski2024}.
The Weiszfeld algorithm solves this efficiently for Euclidean distances~\cite{cooper1981,weiszfeld1937} and the applicability of generalized versions~\cite{eckhardt1980,drezner2009} depend on the structure of $m(\cdot)$.
We thus resort to an iterative approach to solve~\eqref{eq:location_step_problem} by employing a general-purpose optimizer readily available in many programming languages.
For large values of $N_{pc}$, the location step is costly.
However, greedy initialization tends to initialize near local minima in these cases.

In the allocation step, each $\y_i \in W$ is assigned to the preconditioner minimizing Krylov iterations, and we obtain a new partition $W_k$ of $W$.
This defines a \emph{generalized Voronoi diagram}~\cite{chew1985}.

\begin{remark}\label{rem:pcnotatcol}
    Differently from previous work~\cite{graham2021}, we do not assume $\{\hat{\y}_k\}_{k=1}^{N_{pc}} \subset W$.
    Restricting the preconditioners to the parameter locations may be suboptimal and allowing the preconditioners at any point in the parameter enables continuous optimization.
\end{remark}

\subsection{Determining the number of preconditioners}\label{subsec:determining-Npc}
Determining $N_{pc}$ resembles a clustering problem.
Instead, we re-use the greedy initialization: at each step, we estimate the total computation time
\begin{align*}
    \tau_{est}(N_{pc})=\sum_{k=1}^{N_{pc}}\left(N_{ratio}+  \sum_{\y_j\in W_k} m(\y_j - \y_k)\right)
\end{align*}
and, if $\tau_{est}(N_{pc}) > \tau_{est}(N_{pc} - 1) > \tau_{est}(N_{pc} - 2)$, we discard the last two placed preconditioners and set $N_{pc} = N_{pc} - 2$.
Although this procedure has no rigorous guarantee to find the optimal value for $N_{pc}$, it gave us satisfactory results at a cheap cost.
We note that other alternatives are available~\cite{thorndike1953,rousseeuw1987,tibshirani2001}, but these require multiple expensive clustering runs.

We summarize the preconditioner placement in Algorithm~\ref{alg:location-allocation}.

\begin{algorithm}[h]
    \caption{Preconditioner placement}\label{alg:location-allocation}
    \phantom{}\textbf{Input} $W$, $N_{ratio}$, $pc_{fixed}=\{\}$\\
    \phantom{}\Comment{Parameter locations, Cost ratio, Precomputed preconditioner locations}\\
    \phantom{}\textbf{Output} $pc_{loc}$
    \Comment{Computed preconditioner locations}
    \begin{algorithmic}[1]
        \State $N_{pc} \gets |pc_{fixed}| - 1$
        \Comment{The number of preconditioners}
        \State $pc_{loc} \gets pc_{fixed}$
        \Comment{Preconditioner locations}
        \State $cost \gets \left\{\infty, N_{ratio} * N_{pc}+ \sum_{\y_i \in W}m(\y_i, \yh(\yi|pc_{fixed}))\right\}$
        \item[]\Comment{Expected costs of the preconditioning strategy}
        \item[]
        \While{$cost$ not increasing twice in a row}
        \Comment{\textbf{Initialization} and compute $N_{pc}$}
            \State $N_{pc} \gets N_{pc} + 1$
            \State $pc_{loc} \gets pc_{loc} \cup \{\argmax_{\y_i\in W} m(\y_i, \yh(\yi|pc_{loc}) )\}$
            \State $cost\gets cost \cup \{N_{ratio} * N_{pc}+ \sum_{\y_i \in W}m(\y_i, \yh(\yi|pc_{loc}))\}$
        \EndWhile
        \State $N_{pc} \gets N_{pc} - 2, \,\,\,pc_{loc} \gets pc_{loc}[\coloneq2]$
        \Comment{Discard last two preconditioners}
        \item[]
        \While{Preconditioner has changed \AND compute time $<$ strategy gain}
            \item[]\Comment{continue \textbf{location-allocation} as long as there is change and it is worthy}
            \State $W_k \gets partition(W, pc_{loc}, m)$
            \Comment{Compute the partitions $W_k$}
            \For{$\hat{\y}_k \in pc_{loc}$}
                \Comment{Location step}
                \State $\hat{\y}_k \gets \argmin_{{\ys}\in Y} \sum_{\y_i \in W_k[\hat{\y}_k]}m(\y_i-\ys)$
            \EndFor
        \EndWhile\\
        \Return $pc_{loc}$
    \end{algorithmic}
\end{algorithm}

\begin{revenv}
    \begin{remark}[Soft clustering versus hard clustering]
        Contrarily to the hard clustering approach presented here, one could employ soft or fuzzy clustering techniques~\cite{ruspini2019}.
        Although the latter can overcome local minima~\cite{jayaram2013}, they come with higher computational costs~\cite{bora2025} and a hyperparameter that needs to be tuned.
        In Section~\ref{subsec:Npc_experiment}, we will further discuss the effects of local minima in the location-allocation algorithm.
    \end{remark}
\end{revenv}

\begin{remark}[Concentration of measure for isotropic parameter space]\label{rem:high-dimensional-parameter-space}
    As the parameter dimension grows, the distance to the origin of a random vector $X$, and hence the number of Krylov iterations, concentrates~\cite{aggarwal1973}.
    Thus, in high-dimensional isotropic parameter space, the expected number of Krylov iterations with a mean-based preconditioner is nearly uniform.
    If this number is below the cost ratio, mean-based preconditioning is optimal; otherwise, separate preconditioners per parameter value are preferred.
    Consequently, in this case the optimal value of $N_{pc}$ is either one if $m\left(\mathbb{E}[\|X\|_{2}]\right) < N_{ratio}$, or the cardinality of $W$.
\end{remark}

    \section{Estimating $m(\cdot)$ for GMRES}\label{sec:estimating-m}
    To apply Algorithm~\ref{alg:location-allocation}, we need a function approximating the number of Krylov iterations.
We propose an approach based on Gaussian process regression (GPR), using a priori bounds on the number of iterations to initialize a gray-box GPR~\cite{astudillo2022}.
We illustrate this on GMRES iterations with LU preconditioning, setting $\mathbb{P}(\yi)=\mathbb{A}(\yh(\yi))^{-1}$ to solve
\begin{equation}
    \mathbb{A}(\yh(\yi))^{-1}\mathbb{A}(\yi)\u(\yi)=\mathbb{A}(\yh(\yi))^{-1}b,\qquad \text{ for }\yi\in W.\label{eq:precondsystem}
\end{equation}
Next, we construct the gray-box GPR in Sections~\ref{subsec:gray-box-gpr}--\ref{subsec:hyperparameter-tuning} and the training of the GPR is presented in Section~\ref{subsec:gpr-training}.
This surrogate can then be used in the preconditioner placement strategy outlined in Section~\ref{sec:placement}.

\subsection{Gray-box GPR}\label{subsec:gray-box-gpr}
To infer $m(\cdot)$, we use an active learning strategy based on Gaussian process regression and an acquisition function that balances knowledge gain and training cost.
The uncertainty band provided by the Gaussian process, and the GPR will allow us to define a termination criterion.
We consider a Gray-box GPR~\cite{astudillo2022}, modelling $m$ as a nonlinear function of a Gaussian process, with which we include prior knowledge in our estimate.

For GMRES, we use the Elman estimate~\cite{elman1982,graham2021}, which states that the relative tolerance $\varepsilon$ at the $m^\text{th}$ GMRES iteration is bounded by
\begin{align}
    \varepsilon = \frac{\|r_m\|}{\|r_0\|}\leq \left( \frac{2\sqrt{\alpha}}{\alpha + 1} \right)^m\label{eq:alphabound},
\end{align}
with $\|I-\mathbb{A}^{-1}(\yh(\yi))\mathbb{A}(\y_i)\| \leq \alpha < 1$, if we are solving equation~\eqref{eq:precondsystem}.
Therefore, an upper bound on $\|I-\mathbb{A}^{-1}(\yh(\yi))\mathbb{A}(\y_i)\|$ provides us with an upper bound on the required GMRES iterations.
Moreover, if we estimate $\alpha$ as a function of $\y_i-\hat{\y}(\yi)$, equation~\eqref{eq:alphabound} provides us with an estimate for $m(\y_i-\hat{\y}(\yi))$.
Hence, we rewrite~\eqref{eq:alphabound} and define
\begin{equation}
    m(\y)\coloneqq\ln(\varepsilon)\ln\left( \frac{2\sqrt{\mathbb{E}\left[\alpha(\y)\right]}}{\mathbb{E}\left[\alpha(\y)\right] + 1} \right)^{-1}\label{eq:mdef},
\end{equation}
for $\y\in Y$ and where, exploiting translation invariance, $\alpha(\y)$ models $\alpha$ in terms of $\y$, sampled from a Gaussian process
\begin{equation}
    \alpha(\y) \sim GP(\mu_0(\y, C), K(\y, \y'))\label{eq:alpha_GP_def},
\end{equation}
with prior mean $\mu_0(\y, C)$ possibly dependent on $N_{hyp}$ hyperparameters $C=( C_1, \ldots,$ $ C_{N_{hyp}} )\in\mathbb{R}^{N_{hyp}}$, and kernel $K(\y, \y')$.
The prior mean is problem-specific, and needs to be chosen carefully.
In Section~\ref{sec:helmholtz}, we will expand on this in the case of Helmholtz scattering.
Moreover, we elaborate on the kernel choice in Section~\ref{subsec:kernel}.
To simplify the notation, we set
\begin{equation}
    g(\alpha) \coloneqq \ln(\varepsilon)\ln\left(\frac{2\sqrt {\alpha}}{\alpha + 1}\right)^{-1}, \text{\qquad such that \qquad} m(\y) = g\left(\mathbb{E}\left[ \alpha(\y) \right]\right)\label{eq:m_estimate}.
\end{equation}
\begin{remark}[Nonlinearity of $g(\alpha)$]
    Ideally, we would like to use $m(\y) = \mathbb{E}\left[g( \alpha(\y)) \right]$, as this accurately describes the expected value.
    However, this requires integrating over the probability space, which is analytically intractable and computationally expensive.

    To perform the preconditioner placement outlined in Section~\ref{sec:placement} accurately, it is important for $m(\cdot)$ to be accurate near the boundaries of the partition $\left\{ W_k \right\}_{k}$: this determines the distance between the preconditioners, and therefore, $N_{pc}$.
    At the boundaries of the partitions, the values of $\alpha$ are neither very low or high.
    As we can see in Fig.~\ref{fig:gplot},  $g(\alpha)$ is concave for small ($\alpha < 0.01$) values of $\alpha$, and convex for large.
    Outside of these regions, that is in the area $[0.01, 0.03]$ in Fig.~\ref{fig:gplot}, $g(\alpha)$ remains rather linear, such that the effect of switching $\mathbb{E}\left[ \cdot \right]$ and $g(\cdot)$ is minimal.
\end{remark}

\begin{figure}
    \centering
    \begin{tikzpicture}
        \begin{axis}[
        domain=0.001:0.6,
        samples=400,
        xlabel={$\alpha$},
        ylabel={$g(\alpha)$},
        grid=minor,
        thick,
        xmin=0,
        xmax=0.6,
        ymin=0,
        width=0.9\textwidth,
        ylabel near ticks,
        height=0.4\textwidth,
        yticklabels={}, 
        ]
            \addplot[black, thick] ({x}, {-1/ln((2*sqrt(x))/(x+1))});

            \coordinate (spyviewer) at (rel axis cs:0.19,0.98);
        \end{axis}
        \begin{axis}[
            domain=0.00001:0.025,
        at={(spyviewer)},
        anchor={north},
        footnotesize,
        xmin=0,
        xmax=0.025,
        xtick={0, 0.01, 0.02},
        xticklabels={0, 0.01, 0.02}, 
        ymin=0.2,
        ymax=1.01,
        scaled x ticks=false,
        yticklabels={}, 
        ]
            \addplot[black, thick] table [x=X,y=Y,col sep=comma] {data/g_alpha_detail.csv};
        \end{axis}
    \end{tikzpicture}
    \caption{Plot of $g(\alpha)$ from 0 to $0.6$, up to scaling by hyperparameters.}
    \label{fig:gplot}
\end{figure}

\begin{remark}[Different estimates]
    The Elman estimate is based on the field of values of $\mathbb{A}(\yh(\yi))^{-1}\mathbb{A}(\yi)$, but different estimates could be considered, based on eigenvalues or pseudospectra, see~\cite{embree1999} for a clear overview.
    However, different estimates all have their strengths and weaknesses, and we considered the Elman estimate because of the available bounds for the field of values for the Helmholtz equation.
\end{remark}

\begin{revenv}
    \begin{remark}[Application to other PDEs]
        When discussing the estimates outlined in this section, one might wonder how this approach works for different equations. 
        For other PDEs, the overall strategy would be generally applicable.
        However, in Section~\ref{subsec:prior-mean} we do rely on analytical estimates on the number of GMRES iterations~\cite{graham2021}, which are specific to the Helmholtz equation.
        These estimates are key for initializing and thus low training costs for the GP\@.
        To the best of our knowledge, such estimates are not (yet) available in the literature for many other PDEs.
        While one could still use an uninformed prior in these cases, the GP training could potentially be a more expensive strategy.
    \end{remark}
\end{revenv}

\subsection{Kernel}\label{subsec:kernel}
For our method to work efficiently in high parameter dimensions, we have to be careful in choosing the kernel.
The \emph{curse of dimensionality} poses a challenge in Gaussian processes.
Kernel selection may seem like an art, and~\cite[Chapter 2]{duvenaud2014} provides a useful guide in the `Kernel cookbook', while an overview of dimensionality issues can be found in~\cite{binois2022}.
For computational efficiency, we consider univariate kernels~\cite{duvenaud2011,neal1997}:
\begin{equation}
    K(\y, \y') = \sum_{i=1}^N K_i(y_i, y_i')\label{eq:sumkernel}.
\end{equation}
\begin{newenv}
    Each kernel $K_i(\cdot, \cdot)$ is a symmetrized Matérn kernel~\cite{duvenaud2014}:
    \begin{align}
        \nonumber K_i(y_i, y_i') &= K^\text{sym}_i(K^{\text{Mat\'ern}}_{\nu}(y_i,y_i')) \\
        &= \sum_{(m_1,m_2)\in\left\{ -1, 1 \right\}^2}  K^{\text{Mat\'ern}}_{\nu}(m_1 y_i,m_2 y_i'),\label{eq:univariatekernel}
    \end{align}
    where $K^{\text{Mat\'ern}}_{\nu}(\cdot, \cdot)$ is the standard Matérn kernel.
    The symmetrization has been performed using the \emph{sum over orbits} technique~\cite{duvenaud2014}, and to model the nondifferentiability at the origin, we choose $\nu=\frac{1}{2}$ for the smoothness parameter.

    Although this additive approach partly mitigates the high dimensionality, we take another step.
    We take inspiration from~\cite{hvarfner2024} and choose the correlation lengths of the kernels $K_i(\cdot, \cdot)$ to be inversely proportional to the `importance' in the $i^\text{th}$ dimension.
    By increasing the correlation lengths in the parameter dimensions that bear less importance, the variance of the trained Gaussian process will be lower in those dimensions~\cite{shende2022}.
    In essence, this anisotropy mitigates the possible high number of parameters by reducing the effective size of the parameter domain $Y$, reducing the effective model complexity~\cite{hvarfner2024}.
    To achieve this, we define the vector $\{\gamma_i\}_{i=1}^N=\bm{\gamma}\in\mathbb{R}^N$ such that it corresponds to the `importance' of the dimensions.
    A similar approach regarding sparse grids has been taken in~\cite{addy2025}.
\end{newenv}
In Section~\ref{sec:helmholtz}, we outline how we can obtain these anisotropy weights.

\subsection{Hyperparameter tuning}\label{subsec:hyperparameter-tuning}
In the previous sections, we have introduced hyperparameters in the prior mean.
The standard approach of selecting them is by maximizing the likelihood of the training data~\cite{jones1998,rasmussen2000} given the prior mean and kernel.
However, since we use correlation lengths inversely proportional to the importance of the dimension to reduce computational cost, large correlation lengths negatively affect the accuracy of the chosen hyperparameters.

To address this, we choose hyperparameters minimizing the mean squared error of the training data with respect to the prior mean:
\begin{equation*}
C=(C_1, \ldots, C_{N_{hyp}}) = \argmin_{C^\ast \in \mathbb{R}^{N_{hyp}}} \sum_{\y_{i_{train}}}\left(\mu_0(\y_{i_{train}}| C^\ast) - \alpha_{\y_{j_{train}}}  \right)^2,
\end{equation*}
where $\y_{i,{train}}$ are the training locations and $\alpha_{\y_{i,{train}}}$ the corresponding values for $\alpha$.
We solve this minimization problem using a gradient-based optimization routine, as the objective function is smooth and has explicit derivatives.

\subsection{Training the Gray-box GPR}\label{subsec:gpr-training}
In our active learning strategy, we initialize a mean-based preconditioner at $\bar{\y}$, the center of $W$, and solve the linear system~\eqref{eq:linear_system} using the mean-based preconditioner $\mathbb{A}(\bar{\y})^{-1}$ for training points $\y_{i,train} \in W$, recording $\left\{y_{i,train}, \alpha(\y_{i,train})  \right\}$.
These evaluations provide both linear system solutions and training pairs.
Following~\cite{astudillo2022}, we iteratively select training points from $W$ by maximizing an acquisition function $f_{acq}:W\to\mathbb{R}$.

Since evaluation costs depend on $m(\y)$, we account for this in $f_{acq}$, similarly to cost-aware Bayesian optimization~\cite{luong2021,xie2024}.
We define $f_{acq}$ as the variance-to-cost ratio, capped at a cutoff-value $m_{\max}=\frac{\tau_{pc}}{\tau_{GMRES}}$, ensuring unimportant training data is not selected.
The values of $\tau_{pc}$ can be calculated when computing the mean-based preconditioner and $\tau_{GMRES}$ can be computed on the fly using the computed training points.
This is a particular instance of the well-studied expected improvement per unit cost~\cite{snoek2012}:
\begin{equation}
    f_{acq}(\y; m_{max})\coloneqq
    \begin{cases}
        \frac{\mathbb{V}\left[ g(\alpha(\y)) \right]}{\mathbb{E}\left[ g(\alpha(\y)) \right]} &\text{ for } \mathbb{E}\left[ g(\alpha(\y)) \right] \leq m_{max},\\
        -\infty &\text{ else},
    \end{cases}\,\,\,\,\,\,\,
    \label{eq:f-acq-def}
\end{equation}
where $\mathbb{E}\left[ g(\alpha(\y)) \right]$ and $\mathbb{V}\left[ g(\alpha(\y)) \right]$ are the mean and variance of $g(\alpha(\y))$ from~\eqref{eq:m_estimate}, respectively.
To compute~\eqref{eq:f-acq-def}, we approximate $\mathbb{V}\left[ m(\y) \right]$ by central differences:
\begin{equation*}
    \mathbb{V}\left[ g(\alpha(\y)) \right] \approx \frac{g(\mathbb{E}\left[ \alpha(\y) \right] + \mathbb{V}\left[ \alpha(\y) \right])-g(\mathbb{E}\left[ \alpha(\y) \right] - \mathbb{V}\left[ \alpha(\y) \right])}{2}.
\end{equation*}
This will introduce a small bias but enables for efficient computations as the expected values can be computed directly~\cite[Chapter~6]{bishop2006}.

Each new training point is then selected by
\begin{equation*}
    \y_{i,train} = \argmax_{\y_i\in W} f_{acq}(\y_i-\bar{\y};m_{\max}),
\end{equation*}
and each time we add a training point, we fit the hyperparameters again.
The acquisition function requires the hyperparameters itself and hence, to start the training, we initialize the Gaussian process with two points, $\bar{\y}$ and $\y_{\min}=\argmin_{\y\in W} \mu_0(\y, 1)$.
At $\bar{\y}$, we set $\alpha(\bm{0})=g^{-1}(1)$, and at $\y_{\min}$ we solve $\mathbb{A}(\bar{\y})^{-1}\mathbb{A}(\y_{\min})=\mathbb{A}(\bar{\y})^{-1}b$ to obtain $m({\y_{\min}}-\yb)$.

Training continues until the surrogate model reaches a satisfactory accuracy.
Ideally, accuracy is assessed using a test set, but this is computationally expensive.
Instead, we use the \emph{Stabilizing predictions} (SP) stopping criterion~\cite{bloodgood2009,cohen1960}, stopping training when the averaged disagree ratio drops below $1\%$.
Since this criterion is originally defined on labelled data, we adapt it to the continuous scale of $m(\cdot)$ by defining an `agree' whenever the relative difference between two iterations is lower than $1\%$, with a minimum of 1 iteration.

\subsection{Assembling the surrogate}\label{subsec:assembling-the-surrogate}
With the collected training points $W_{train}=\left\{ \y_{1,train}, \ldots, \y_{N_{train}, train} \right\}$ and obtained Gaussian process values $\bm{\alpha}_{train}$$=\{ \alpha_{\y_{1, train}} ,$ $ \cdots,\alpha_{\y_{N_{train}, train}} \}$, we define $m(\cdot)$ as~\eqref{eq:mdef}, with
\begin{equation*}
    \mathbb{E}\left[\alpha(\y)\right] = \mu_0(\y) + K(\y+\yb, W_{train})^\top K(W_{train},W_{train})^{-1}(\bm{\alpha}_{train}-\mu_0(W_{train})),
\end{equation*}
where $\mu_0(\cdot)$ is the prior mean from~\eqref{eq:GP-mean-def_y} and $K(\cdot, \cdot)$ is the kernel~\eqref{eq:sumkernel}.
We summarize our strategy in Algorithm~\ref{alg:gray-box-gpr}.

\begin{algorithm}
    \caption{Gray-box GPR $m(\y)$ training}\label{alg:gray-box-gpr}
    \phantom{}\textbf{Input} $W$, $\bar{\y}$
    \Comment{Parameter locations, Mean-based parameter location}\\
    \phantom{}\textbf{Output}
    \phantom{}\hspace{\algorithmicindent} $Y_{eval}$, $m(\cdot)$, $m_{\max}$\\
    \phantom{}\Comment{Evaluated locations, $m(\cdot)$ function, maximal valid value of $m(\cdot)$}
    \begin{algorithmic}[1]
        \State $\alpha \gets GP(\mu_0(\cdot), K(\cdot, \cdot)(\cdot, \cdot))$
        \Comment{Initiate GP~\eqref{eq:alpha_GP_def}}
        \State $\alpha \gets (\bm{0}, g^{-1}\left(1\text{ iteration}\right))$
        \State $\tau_{PC} \gets PC(\bar{\y})$
        \Comment{Initialize mean-based PC and store computation time}
        \State $\y_{\min} = \argmin_{\y\in W} \mu_0(\y, 1)$
        \State $Y_{eval} \gets \{\y_{\min}\}$
        \Comment{List containing all evaluated parameter locations}
        \State $\tau_{\y_{\min}}, m_{\y_{\min}} \gets GMRES(\mathbb{A}(\bar{\y})^{-1}\mathbb{A}(\y_{\min})=\mathbb{A}(\bar{\y})^{-1}b)$
        \Comment{Initialize by evaluating for smallest $\|\cdot\|_{\gamma, 2}$}
        \State $\tau_{tot}, m_{tot} \gets \tau_{\y_{\min}}, m_{\y_{\min}}$
        \Comment{Total GMRES iterations computed and time spend on it}
        \State $\alpha \gets (\y_{\min}-\bar{\y}, g^{-1}\left(m_{\y_{\min}}\right))$
        \Comment{Add training data to the GP}
        \State $SP \gets \{\}$
        \Comment{List for SP results}
        \item[]
        \While{running average($SP$) $>$ ($1\%$ and 1)}
            \State $m_{\max} \gets {\tau_{PC}}/{\left( \tau_{tot}  / m_{tot}\right)}$
            \Comment{Update the maximal relevant value of $m(\cdot)$}
            \State $\y_{train} \gets \argmax\left\{ f_{acq}(\y_i-\bar{\y}; m_{\max}): \y_i \in (W\setminus Y_{eval}),  \right\}$
            \Comment{New train point}
            \State $Y_{eval} \gets Y_{eval} \cup \{\y_{train}\}$
            \State $\tau_{\y_{train}}, m_{\y_{train}} \gets GMRES(\mathbb{A}(\bar{\y})^{-1}\mathbb{A}(\y_{train})=\mathbb{A}(\bar{\y})^{-1}b)$
            \State $\alpha \gets (\y_{train}-\bar{\y}, g^{-1}\left(m_{\y_{train}}\right))$
            \Comment{Add training data to GP}
            \State $SP \gets SP \cup \{SP(g\left(\mathbb{E}\left[  \alpha  \right]\right))\}$
            \Comment{Compute disagreement ratio}
            \State $\tau_{tot}\gets \tau_{tot}+\tau_{\y_{\min}}$
            \Comment{Bookkeeping}
            \State $iter_{tot} \gets m_{tot}+m_{\y_{\min}}$
        \EndWhile
        \item[]\\
        \Return $Y_{eval}, g\left(\mathbb{E}\left[  \alpha(\cdot)  \right]\right), {\tau_{PC}}/{\left(\tau_{tot}  / m_{tot}\right)}$
    \end{algorithmic}
\end{algorithm}

With this approximation $m(\cdot)$ satisfying Assumption~\ref{ass:shift_invariant}, we can use Algorithm~\ref{alg:location-allocation} to locate the preconditioner locations for the remaining parameter locations.
Since a mean-based preconditioner was placed during the Gray-box GPR training, we pass this information on to the location-allocation procedure.
Specifically, we run Algorithm~\ref{alg:location-allocation} with $W=W\setminus\left\{ Y_{eval} \right\}$, $N_{ratio}=m_{\max}$ to discard the already evaluated preconditioners, and the mean-based preconditioner $pc_{fixed}=\{\bar{\y}\}$.

\begin{remark}[Parallel training]
    Instead of our sequential training approach, parallel approaches could be explored.
    See, for example,~\cite{ginsbourger2010,snoek2012}.
\end{remark}

    \section{Helmholtz problem}\label{sec:helmholtz}
    We study the Truncated Exterior Dirichlet Problem (TEDP) of the Helmholtz equation in two dimensions, where an incoming wave scatters from an impenetrable scatterer with boundary $\Gamma_{in}$.
We truncate \begin{revenv} the \end{revenv} unbounded domain with a circle $\Gamma_{out}$ of radius $r_{out}$ enclosing the scatterer.
We use a Robin boundary condition~\cite{shirron1998} on the scattered wave to approximate the \emph{Dirichlet-to-Neumann map}~\cite[Section 2.6.3]{nedelec2001}.
Denoting $D=\left\{ \x\in\mathbb{R}^2\,:\,r_{in}(\theta(\x)) < \|\x\| < r_{out} \right\}$, we thus consider the variational problem: for $\y\in Y$, find $u(\y)\in H^1(D)$ such that
\begin{align}
    \int_D A(\y) &\nabla u(\y) \cdot \nabla v \dx - k_0^2\int_{D}n(\y) u(\y) v \dx - i k_0\int_{\Gamma_{out}} u(\y) v\dx \nonumber \\ &=\int_{\Gamma_{out}} \left( \frac{\partial}{\partial \hat{\bm{n}}} -ik_0\right)u_{in} v\dx, \qquad \qquad
    \text{for all }v \in H^{1}(D)\label{eq:weakform},
\end{align}
where $u_{in}$ is the incoming wave, $n(\y; \cdot)$ is a heterogeneous, real-valued refractive index of the medium, and $A(\y; \cdot)$ a heterogeneous coefficient, taking values in $\mathbb{R}^{2\times 2}$.
Moreover, we require both $n(\y,\cdot)$ and $A(\y,\cdot)$ to be $C^1$ in $\y$ to go beyond affine expansions in $\y$.

Discretization with the finite element method leads to linear systems
\begin{equation}
    \mathbb{A}(\y)\u(\y)=b.\label{eq:linear_system}
\end{equation}
To use a-priori bounds from~\cite{graham2021}, we assume that problem~\eqref{eq:weakform} is nontrapping, where we denote $\left\| \cdot \right\|_{H_{k_0}(D)}^2\coloneqq \left\| \nabla \cdot \right\|_{L^2(D)}^2+k_0^2\left\| \cdot\right\|_{L^2(D)}^2$:
\begin{assumption}[Nontrapping]\label{ass:nontrapping}
The quantities $D$, $A$, and $n$ are such that, given $f\in L^2(D)$, the solution $u(\y)$ of~\eqref{eq:weakform} exists, is unique, and satisfies
\begin{align*}
    \|u(\y)\|_{H_{k_0}^1(D)}\leq C_{\text{bound}} \|f\|_{L^2(D)},
\end{align*}
where $C_{bound}$ is independent of $k_0$ and $\y\in Y$.
\end{assumption}

For high $k_0$, solving~\eqref{eq:linear_system} becomes challenging with GMRES\@.
To solve~\eqref{eq:linear_system} many times for different parameter values, we follow Section~\ref{sec:placement} to assign the linear systems to a limited number of preconditioners using the surrogate outlined in Section~\ref{sec:estimating-m}.
Therefore, we need a good approximation for the prior mean, which we will construct in the next section, and a measure of dimension importance, which we will discuss in Section~\ref{subsec:importance-weights}.

\subsection{The prior mean}\label{subsec:prior-mean}
Choosing a prior mean is a delicate task.
We could use the estimates in~\cite{graham2021} and set
\begin{equation}
    \mu_0(\y,C)\coloneqq C_1 \|A(\y; \cdot) - A(\hat{\y};\cdot)\|_{L^2_{op}(D)}+C_2 \|n(\y; \cdot) - n(\hat{\y};\cdot)\|_{L^2(D)}\label{eq:GP-mean-def},
\end{equation}
where $C_1$ and $C_2$ are hyperparameters and $\|\cdot\|_{L^2_{op}(D)}$ is the $L^2$-norm of the pointwise spectral 2-norm.
However, computing this requires integration over $D$, which is costly.

To evaluate~\eqref{eq:GP-mean-def} efficiently, we exploit the differentiability of $A$ and $n$:
\begin{equation*}
    n(\y, \cdot) = n(\yh, \cdot) + \nabla_{\y} n(\bm{\eta}^1, \cdot) \cdot \left( \y - \yh \right),
\end{equation*}
for some $\bm{\eta}^1\in Y$ between $\y$ and $\yh$.
Expanding the second term in~\eqref{eq:GP-mean-def}, we obtain:
\begin{align*}
    \|n(\y; \cdot) - n(\hat{\y};\cdot)\|_{L^2(D)}&= \| \nabla_{\y} n(\bm{\eta}^1, \cdot) \cdot \left( \y - \yh \right) \|_{L^2(D)}\\
    &= \left( \left(\y - \yh \right)^\top\mathbb{B}(\bm{\eta}^1)\left(\y - \yh \right) \right)^{\frac{1}{2}}
    = \left\| \y-\yh \right\|_{\mathbb{B}(\bm{\eta}^1)},
\end{align*}
where $\left\| \cdot \right\|_{\mathbb{B}(\bm{\eta}^1)}$ denotes a \emph{weighted $\ell^2$-norm} with a weight matrix $\mathbb{B}(\bm{\eta}^1)$, parameterized by $\bm{\eta}^1$.
In this case, $\mathbb{B}(\bm{\eta}^1)$ is a symmetric matrix with entries
\begin{equation}
    \mathbb{B}_{ij}(\bm{\eta}^1) = \int_D \partial^{e_i} n(\bm{\eta}^1, \cdot)\partial^{e_j} n(\bm{\eta}^1, \cdot)  \dx \label{eq:Bmatdef}.
\end{equation}
Similarly, we obtain $\|A(\y; \cdot) - A(\hat{\y};\cdot)\|_{L^2_{op}(D)} = \left\| \y-\yh \right\|_{\mathbb{D}(\bm{\eta}^2)},$ where $\mathbb{D}(\bm{\eta}^2)$ is a symmetric matrix with entries
\begin{equation}
    \mathbb{D}_{ij}(\bm{\eta}^2) = \int_D\left\|  \partial^{e_i} A(\bm{\eta}^2, \cdot)\partial^{e_j} A(\bm{\eta}{^2}, \cdot)\right\|_{2,2}  \dx \label{eq:Dmatdef},
\end{equation}
where $\bm{\eta}^2\in Y$ lies between $\y$ and $\yh$ and $\left\| \cdot \right\|_{2,2}$ is the spectral matrix norm.

The matrices $\mathbb{B}(\cdot)$ and $\mathbb{D}(\cdot)$ can often be approximated by constant matrices by exploiting the structure of the parameter dependence.
We will discuss this for an affine expansion in Section~\ref{subsec:disjoint-density}, and for a parameterized domain in Section~\ref{subsec:parameterized-domain}.
However, the results presented in both sections are more widely applicable than these two examples, which merely show the process of applying these methods.
We also note that, because of the hyperparameters $C_1$ and $C_2$, it is sufficient to know them up to an (unknown) constant.

We can then combine~\eqref{eq:GP-mean-def} with $\mathbb{B}$ and $\mathbb{D}$ to obtain:
\begin{equation}
    \mu_0(\y|C)= C_1 \left\| \y-\yh \right\|_{\mathbb{D}}+C_2  \left\| \y-\yh \right\|_{\mathbb{B}}\label{eq:GP-mean-def_y},
\end{equation}
for the prior mean of the gray-box Gaussian process~\eqref{eq:alpha_GP_def}.

\subsection{Anisotropy weights}\label{subsec:importance-weights}
Next to the prior mean, we need a value for the anisotropy weights $\gamma_j$ proportional to the importance of dimensions in the parameter space.
Therefore, we define them in terms of the matrices $\mathbb{B}$ and $\mathbb{D}$:
\begin{equation*}
    \gamma_j \coloneqq C_1\sqrt{\mathbb{D}_{ii}}+C_2\sqrt{\mathbb{B}_{ii}},
\end{equation*}
where $C_1$ and $C_2$ are the hyperparameters from equation~\eqref{eq:GP-mean-def}.
Now, the correlation length $l_i$ of the Matérn kernel in dimension $i$ is given by:
\begin{equation*}
    l_i = \text{diam}(D) \frac{\max_i \gamma_i}{\gamma_i} \geq \text{diam}(D),
\end{equation*}
where $\text{diam}(D)$ is the diameter of the computational domain.

To validate our methods, we will work through two applications.
In the next section, Section~\ref{subsec:disjoint-density}, we will consider an affine dependence on the parameter $\y$ for $n$ and in Section~\ref{subsec:parameterized-domain}, we consider a parameterized domain test case, where shape of the central scatterer depends on the parameter $\y$.

\subsection{Affine expansion}\label{subsec:disjoint-density}
We investigate affine dependence on the parameter $\y$ for $n$, and a constant $A=I_2$.
The affine expansion of the refractive index $n$ is given by
\begin{equation}
    n(\y,\x) \coloneqq n_{0}(\x) + \sum_{i=1}^{N} \psi_i(\x) y_i, \label{eq:affinendef}
\end{equation}
with $n(\y, \x)$ such that the resulting problem is a non-trapping problem for all $\y\in Y$ and $\left\{ \psi_i \right\}$ is uniformly positive over $\x$.

To apply the methods we developed in Section~\ref{sec:estimating-m}, we need to estimate the matrices $\mathbb{D}$ and $\mathbb{B}$.
Since $A$ is independent of $\y$, we set $\mathbb{D}_{ij}=0$.
To estimate $\mathbb{B}$, we use equation~\eqref{eq:Bmatdef} and compute
\begin{equation}
    \mathbb{B}_{ij} =\int_D \psi_i(\x)\psi_j(\x) \dx.\label{eq:affineBdef}
\end{equation}
In the numerical experiments, we will consider the expansion
\begin{equation}
    n(\y,\x) \coloneqq 1 + \sum_{i=1}^{N} \mathbbm{1}_{\Omega_i}(\x)\chi(\x)\eta_i \frac{y_i-1}{2}, \label{eq:ndef}
\end{equation}
for a partition $\{\Omega_i\}_{i=1}^N$ of $D$, weights $\eta_i\in\mathbb{R}$, $\y\in Y=[-1,1]^N$, and $\chi$ a mollifier satisfying the assumption below:
\begin{assumption}\label{ass:mollifier}
The mollifier $\chi(\x)$ is continuous on $D$, $\chi(\x)=1$ on $\Gamma_{in}$ and $\chi(\x)=0$ on an open neighbourhood of $\Gamma_{out}$.
\end{assumption}

By employing a mollifier, we ensure that $n(\cdot, \x)\equiv 1$ and $\nabla_{\x} n(\y,\x)=0$ on the outer boundary $\Gamma_{out}$, to obey the premise of the Sommerfeld radiation condition.
Moreover, we choose the partition $\{\Omega_i\}_{i=1}^N$ of $D$ as
\begin{equation*}
    \Omega_i = \left\{\x \in D \,,\, \frac{2\pi i}{N} \leq \theta(\x) < \frac{2\pi (i+1)}{N} \right\},
\end{equation*}
such that, together with the term $(y_i-1)$ in~\eqref{eq:ndef}, the non-trapping Assumption~\ref{ass:nontrapping} is fulfilled~\cite[Condition~2.6]{graham2019}.
Moreover, $\mathbb{B}$ is a diagonal matrix with entries $\mathbb{B}_{ii}=C_\chi\eta_i^2$
where $C_\chi \in \mathbb{R}$ is a constant depending on the mollifier that will be absorbed into the hyperparameters.
This reduces the correlation lengths to
\begin{equation*}
    l_i = \diam(D) \frac{\max_j \eta_j}{\eta_i}.
\end{equation*}
For a generic $n$, the matrix $\mathbb{B}$ from~\eqref{eq:affineBdef} can be computed similarly to the techniques we present in the next section.

\subsection{Parameterized domain}\label{subsec:parameterized-domain}
Our second test case considers scattering against a star-shaped object $\mathcal{D}_{scat}(\y)$, parameterized by $\y\in Y = \left[ -1,1 \right]^N$.
Due to this star-shaped property, we can express $\mathcal{D}_{scat}(\y)$ in polar coordinates as the region inside a simple, closed curve $r(\y, \theta)$, $\theta \in [0, 2\pi]$:
\begin{equation*}
    \mathcal{D}_{scat}(\y) = \left\{ \x \in \mathbb{R}^2 : |\x| \leq r(\y, \theta(\x)) \right\},
\end{equation*}
where $r(\cdot, \theta)\in C^{0,1}_{per}([0, 2\pi))$ for all $\y \in Y$ and $\theta(\x)$ is the angle of $\x$ in polar coordinates.
The resulting scattering problem is posed on the domain
\begin{equation*}
    \mathcal{D}(\y) = \left\{ \x \in \mathbb{R}^2 : r(\y, \theta(\x)) \leq |\x| \leq r_{out} \right\}.
\end{equation*}
We model the boundary $r(\y, \theta)$ as an affine combination of a nominal radius $r_{in}(\theta) \in C^{0,1}_{per}([0, 2\pi))$ and a set of basis functions $\left\{ \psi_j(\theta) \right\}_{j=1}^N \subset C^{0,1}_{per}([0, 2\pi))$:
\begin{equation}
    r(\y, \theta) = r_{in}(\theta) + \sum_{j=1}^N \psi_j(\theta)y_j, \qquad\qquad\text{ for all }\theta\in [0, 2\pi) \text{ and } \y \in Y. \label{eq:affinerdef}
\end{equation}
\begin{revenv}
    This way, the $L^{\infty}$-norm of $\left\{ \psi_j \right\}_{j=1}^N$ encodes the parameter anisotropy.
\end{revenv}

Moreover, we require $\left\{ \psi_j \right\}_{j=1}^N$ to be such that $\min_{\y\in Y}\min_{\theta \in [0, 2\pi)} r(\y,\theta) > 0$ ensuring~\eqref{eq:affinerdef} defines a non-intersecting curve for all $\y \in  Y$.
First, we will consider the formulation on the parameterized domain $\mathcal{D}(\y)$, and then we will reconnect it to the model problem~\eqref{eq:weakform} through a mapping approach.

Similarly to the beginning of Section~\ref{sec:helmholtz}, we use a Robin boundary condition and obtain the variational problem: for all $\y\in Y$, find $\mathrm{u}(\y) \in H^1(\mathcal{D}(\y))$ such that:
\begin{align}
    \int_{\mathcal{D}(\y)}&  \nabla \mathrm{u}(\y) \cdot \nabla v \dx - k_0^2\int_{\mathcal{D}(\y)} \mathrm{u}(\y) v \dx - i k_0\int_{\Gamma_{out}} \mathrm{u}(\y) v\dx \nonumber\\&= \int_{\Gamma_{out}} \left( \frac{\partial}{\partial \hat{n}} -ik_0\right)u_{in} v\dx, \text{\qquad for all }v \in H^{1}(\mathcal{D}(\y))\label{eq:weakformshape},
\end{align}
where, without loss of generality, the $A$ and $n$ from equation~\ref{eq:weakform} are taken constant and set to one such that they define a nontrapping problem~\cite[Theorem~2.5]{graham2019}\footnote{To use this result, we require the boundary excitation to be zero. This is not true for our total field equations, but we can use a bound on the linear combination of the scattered field and the total field such as in~\cite[Corollary~1.6]{chaumont-frelet2023} to obtain the required results.}.
We take care to distinguish $\mathrm{u}$, the solution on the parameterized domain, from $u$, which will later be the solution on a fixed domain.
Together with the star-shaped property of $\mathcal{D}_{scat}(\y)$, we satisfy Assumption~\ref{ass:nontrapping}~\cite[Prop. 3.1, Chapter 5]{lax1989} with the constant $C_{\text{bound}}$ independent of $\y$~\cite{chandler-wilde2008}.

\begin{figure}
    \centering
    \begin{tikzpicture}[thick, scale=1.6]
        \draw[blue,domain=0:3*pi,samples=500,shift={(3,0)}] plot ({0.5*(1.1*cos(\x r)*(1+0.05*cos(3*\x r)+0.04*cos(8*\x r)+0.015*cos(11*\x r))}, {0.5*(1.1*sin(\x r)*(1+0.05*cos(3*\x r)+0.04*cos(8*\x r)+0.015*cos(11*\x r))});
        \draw[green,domain=0:3*pi,samples=500,shift={(3,0)}] plot ({0.5*(1.1*cos(\x r)*(1+0.05*sin(3*\x r)+0.02*cos(6*\x r)+0.01*sin(10*\x r)))}, {0.5*(1.1*sin(\x r)*(1+0.05*sin(3*\x r)+0.02*cos(6*\x r)+0.01*sin(10*\x r)))});
        \draw[red,domain=0:3*pi,samples=500,shift={(3,0)}] plot ({0.5*(1.1*cos(\x r)*(1+0.15*sin(2*\x r)+0.02*sin(6*\x r)+0.01*cos(9*\x r)))}, {0.5*(1.1*sin(\x r)*(1+0.12*sin(2*\x r)+0.02*sin(6*\x r)+0.01*cos(9*\x r)))});
        \draw[orange,domain=0:3*pi,samples=500,shift={(3,0)}] plot ({0.5*(1.1*cos(\x r)*(1+0.15*cos(\x r)+0.01*cos(5*\x r)+0.02*sin(11*\x r)))}, {0.5*(1.1*sin(\x r)*(1+0.15*cos(\x r)+0.01*cos(5*\x r)+0.02*sin(11*\x r)))});
        \draw[fill=none,domain=0:2*pi,samples=500,fill opacity=0.2,shift={(3,0)},dashed] plot ({deg(\x)}:{1.5});

        \draw[->, thick] (-1,0.5) to [bend left=45]  node [above,pos=0.5] {$\left\{ \Phi(\y_i, \cdot) \right\}_{i\in I}$} (1.8,0.5);
        \draw[<-, thick] (-1,-0.5) to [bend right=45]  node [below,pos=0.5] {$\left\{ \Phi^{-1}(\y_i, \cdot) \right\}_{i\in I}$} (1.8,-0.5);
        \node at (-2.15,1) {${D}$};
        \node at (2.85,1) {$\left\{ D(\y_i) \right\}_{i\in I}$};
        \draw[fill=none,domain=0:2*pi,samples=500,fill opacity=0.2,shift={(-2,0)},dashed] plot ({deg(\x)}:{1.5});
        \draw[fill=none,domain=0:3*pi,samples=500,fill opacity=0.2,shift={(-2,0)}] plot ({deg(\x)}:{0.5});
        \draw[-to,shift={(-2,0)}] (0,0) -- (0.703/2, 0.703/2) node [pos=1, right] {$r_0(\theta)$};
    \end{tikzpicture}
    \caption{Schematic overview of the mapping approach. The domain on the left is the reference domain, with the parameter-dependent mapping $\Phi(\y, \cdot)$ from the reference domain to several values of the parameter $\y_i$. The dashed circle has radius $r_{out}$ and signifies the outer radius of the compuatational domain.}
    \label{fig:mappingappraochtikz}
\end{figure}

We use a mapping approach~\cite{tartakovsky2006,xiu2006}, in the spirit of~\cite{castrillon-candas2016,harbrecht2016,hiptmair2018}, on the reference domain to pull the variational formulation~\eqref{eq:weakformshape} back to the reference domain $D=\mathcal{D}(\bm{0})$.
\begin{revenv}
    This approach is shown schematically in Fig.~\ref{fig:mappingappraochtikz} and has the advantage that all different domains $\mathcal{D}(\y)$ can be resolved using the same mesh, which we create for $\y=\bm{0}$ only.
\end{revenv}
We perform the pullback using a parameter-dependent diffeomorphism $\Phi(\y):D\mapsto \mathcal{D}(\y)$, which we consider to be affinely dependent on the entries of $\y$~\cite{cohen2018}:
\begin{equation}
    \label{eq:affinemapping}
    \Phi(\y,\x) = \x + \sum_{j=1}^N \Phi_j(\x)y_j, \qquad \x\in D,
\end{equation}
where $\Phi_j(\cdot) \in W^{1,\infty}(D)$ are the partial transformation maps such that $\Phi_j^{-1}(\y, \cdot)$, $\Phi^{-1}(\y, \cdot) \in W^{1,\infty}(D)$.
With this, the Courant-Fisher theorem for singular values~\cite[Thm.~3.1.2]{horn1991} applied to the Jacobian matrices $\D\Phi(\y)$ and $\D\Phi^{-1}(\y)$ ensures that:
\begin{equation*}
    \sigma_{min} \leq \sigma_1(\y; \x), \cdots, \sigma_d(\y; \x) \leq \sigma_{max}, \quad  \text{for a.e. }\x\in \mathcal{D}(\y)\text{ and for all }\y \in Y.
\end{equation*}
Denoting by $u$ the pulled back solution, we have the variational problem: for every $\y\in Y$, find $u(\y) \in H^1(D)$ such that
\begin{align*}
    \int_{D}& A(\y,\cdot) \nabla u(\y) \cdot \nabla v \dx - k_0^2\int_{D} n(\y,\cdot) u(\y) v \dx - i k_0\int_{\Gamma_{out}}  u(\y) v\dx \\
    &= \int_{\Gamma_{out}} \left( \frac{\partial}{\partial \hat{\bm{n}}} -ik_0\right)u_{in} v\dx +  \int_{D} f v\dx, \nonumber
    \text{\qquad for all }v \in H^{1}(D),
\end{align*}
where the coefficients are $A(\y; \x) =  \D\Phi^{-1}(\y; \x) \D\Phi^{-\top}(\y; \x)\det(\D\Phi(\y; \x))$ and $n(\y; \x) = \det(\D\Phi(\y; \x))$.
Hence, we are considering a parameterized Helmholtz equation of type~\eqref{eq:weakform}, which we solve using a finite element approach to get linear systems of the form ~\eqref{eq:linear_system}.

Similar to the affine expansion in Section~\ref{subsec:disjoint-density}, we still have to obtain the matrices $\mathbb{B}$ and $\mathbb{D}$ by approximating $\eqref{eq:Bmatdef}$ and $\eqref{eq:Dmatdef}$.
To do this, we will use upper bounds in~\cite{harbrecht2016} on the first derivative of $A$ and $n$ with respect to $\y$.

We treat the $\mathbb{B}$ matrix first, by employing~\cite[Lemma~4]{harbrecht2016} to obtain:
\begin{equation*}
    \left|  \partial^{e_i}n(\bm{\eta}^1, \cdot)\right| \leq 2 (1+\sigma_{max})^2 \left\| \Phi_i \right\|_{W^{1,\infty}(D)}, \qquad i=1,\cdots,N,
\end{equation*}
for all $\bm{\eta}^1$, where $\sigma_{\max}$ is the upper bound from the Courant-Fisher Theorem, and $\Phi_i$ are the partial transformations from equation~\eqref{eq:affinemapping}.
Hence, we bound:
\begin{equation*}
    \mathbb{B}_{ij}(\bm{\eta}^1) = \int_D \partial^{e_i} n(\bm{\eta}^1, \cdot)\partial^{e_j} n(\bm{\eta}^1, \cdot)  \dx
    \leq C_{\mathbb{B}} \left\| \Phi_i \right\|_{W^{1,\infty}(D)}  \left\| \Phi_j \right\|_{W^{1,\infty}(D)}
    \eqqcolon C_{\mathbb{B}} \mathbb{B}_{ij},
\end{equation*}
where the constant $C_\mathbb{B}=|D|4 (1+\sigma_{\max})^4$ will be absorbed into the hyperparameter $C_2$.
Similarly, we use ~\cite[Theorem~4]{harbrecht2016} to bound~\eqref{eq:Dmatdef}:
\begin{equation*}
    \mathbb{D}_{ij}(\bm{\eta}^2) \leq C_{\mathbb{D}} \mathbb{D}_{ij},
\end{equation*}
where $\mathbb{D}_{ij}=\mathbb{B}_{ij}=\left\|\Phi_i \right\|_{W^{1,\infty}(D)}  \left\| \Phi_j \right\|_{W^{1,\infty}(D)}$ and $C_\mathbb{D}=\frac{16(1+\sigma_{\max})^2 }{\sigma_{\min}^2} \frac{4 (1 + c_\gamma)^2}{\sigma_{\min}^4\ln\left( 2 \right)^2}$ is a constant which will be absorbed into the hyperparameter $C_1$.
Moreover, $c_\gamma=\sum_{k=1}^N \left\| \Phi_k \right\|_{W^{1,\infty}(D)}$.
In the numerical experiments, we will consider a Fourier expansion for $\left\{\psi_j \right\}_{j=1}^N$:
\begin{equation*}
    \psi_j(\theta)\coloneqq\begin{cases}
                               \vartheta & \text{for }j=1,\\
                               \vartheta \left(  \frac{j+2}{2}\right)^{-\alpha}\sin(j \theta/2 ) & \text{for $j$ even },\\
                               \vartheta \left( \frac{j+1}{2} \right)^{-\alpha}\cos((j-1) \theta / 2 ) & \text{else},\\
    \end{cases}
\end{equation*}
where $\vartheta>0$ is a scaling constant determining the magnitude of the shape deformations and $\alpha$ a parameter that describes the decay in Fourier modes, introducing anisotropy into the parameter domain.
\begin{revenv}
    By choosing this, the lower order harmonics are more influential on the magnitude of the shape variations.
    On the contrary, the higher harmonics have a lower impact on the resulting shape deformation.
\end{revenv}

To make sure that $r(\y, \theta)$ in equation~\eqref{eq:affinerdef} defines a non-intersecting curve for all parameter dimensions $N$, we take
\begin{equation*}
    \vartheta < \vartheta_{\max}(\alpha) \coloneqq \frac{\inf_{\theta\in [0,2\pi)}r_{in}(\theta)}{1 + \sqrt{2}\left(\zeta(\alpha) - 1  \right)},
\end{equation*}
where $\zeta(\alpha)$ is the Riemann Zeta function evaluated at $\alpha$.
Now, combined with the Fourier expansion, we define the partial transformations $\Phi_j(\x)$:
\begin{align*}
    \Phi_j(\x) \coloneqq \chi(\x) \psi_j(\theta(\x)),
\end{align*}
with a mollifier satisfying Assumption~\ref{ass:mollifier} such that we can compute:
\begin{equation}
    \left\| \Phi_j \right\|_{W^{1,\infty}(D)} \leq
    \begin{cases}
        2\vartheta\left\| \nabla \chi(\x) \right\|_{L^{\infty}} & \text{for }j=1,\\
        \left( \frac{j+2}{2} \right)^{-\alpha}\vartheta\left( 1 + \left\| \nabla \chi(\x) \right\|_{L^{\infty}} + \left( \frac{j}{2} \right)  \right)  & \text{for even }j,\\
        \left( \frac{j+1}{2} \right)^{-\alpha}\vartheta\left( 1 + \left\| \nabla \chi(\x) \right\|_{L^{\infty}} + \left( \frac{j-1}{2} \right) \right) & \text{else}.\\
    \end{cases}\label{eq:shapew1infnormdef}
\end{equation}
Finally, we can use the estimates of $\mathbb{B}$ and $\mathbb{D}$ to compute the correlation lengths:
\begin{equation*}
    l_i = \diam(D) \frac{\max_j \left\| \Phi_j \right\|_{W^{1,\infty}(D)}}{\left\| \Phi_i \right\|_{W^{1,\infty}(D)}}.
\end{equation*}

    \section{Numerical illustrations}\label{sec:numerical-investigation}
    To show the efficacy of our method, we consider the two test cases outlined in Section~\ref{sec:helmholtz}.
In these, we set \begin{newenv}
    $r_{in}\equiv 0.4$,
\end{newenv} $r_{out}=1$, and $r_{mol}=0.9 < r_{out}$ to define the mollifier
\begin{equation*}
    \chi(\x) \coloneqq \frac{|\x| - r_{mol}}{r_{in}-r_{mol}},
\end{equation*}
satisfying Assumption~\ref{ass:mollifier}.
We thus have $\diam(D)=2$ and $\left\| \nabla \chi(\x) \right\|_{L^{\infty}}=\frac{1}{r_{mol} - r_{in}}$.
We initialize with a mean-based preconditioner at $\yb=\bm{0}\in\mathbb{R}^N$.
To test our routine, we consider $|W|=1000$ points, uniformly spaced in $Y=[-1,1]^N$.

We implement our numerical methods in Python~3.10, running in an Ubuntu~$24.04$ Docker container.
To obtain the linear systems, we mesh our computational domain using GMesh~\cite{geuzaine2009} with uniform mesh size $h\sim k_0^{-\frac{3}{2}}$ to counteract the pollution effect~\cite{babuska1997}.
We use DOLFINx~\cite{alnaes2014a,baratta2023,scroggs2022,scroggs2022a} to construct the linear systems, which we solve with GMRES in PETSc~\cite{brown2022} until the default relative error ($10^{-5}$).
Moreover, we construct the required preconditioners using PETSc as well.
GMRES iterations and PETSc preconditioner construction have been timed using the \texttt{time.process\_time()} method of the \texttt{time} module in Python.
We perform the optimization during the location-allocation using the L-BFGS-B Fortran subroutine.
When timing pieces of code, we exclude any overhead caused by the assembly of the linear systems~\eqref{eq:linear_system} because we need to assemble them anyway, independently of the preconditioning strategy.
We will first focus on the parameterized shape problem, as these results give more insight into the efficacy of the algorithm.
Subsequently, we will use the affine case to discuss the effects of isotropy in the parameter domain and to see a practical example of the concentration effect.

\subsection{Algorithm performance for the parameterized domain}\label{subsec:conventional_comparison}
We demonstrate the improvements of distributed preconditioning over more traditional approaches like mean-based preconditioning, using parameter dimensions ranging from $N=2$ to $N=25$.

For each experiment, we report the time to train the gray-box GPR ($t_{train}$), the location-allocation time ($t_{l-al}$), and the execution time ($t_{exec}$) of solving the preconditioned linear systems, rounded to the nearest second, and averaged over several runs to reduce the variance.
Moreover, we list the average number of GMRES iterations in the executed strategy ($it_{av}$).
The total computation time ($t_{tot}=t_{train} + t_{l-al}+t_{exec}$) follows.
Additionally, we list the average number of preconditioners $N_{pc}$ over the runs.

\begin{revenv}
    First, we will discuss the algorithm performance with respect to the parameter dimension in Section~\ref{subsec:algperf}, after which we will investigate the effect of $\vartheta$ on the performance.
    In Section~\ref{subsec:anisotropytest}, we test the effect of the parameter anisotropy parameter $\alpha$ on the performance and we conclude this section with a discussion on the effect of the frequency $k_0$ on the performance of the algorithm.

    For a comparison with a different algorithm, we refer the reader to the Section~SM4 in the supplementary material.
    There, we compare our approach with the $K$-means approach as outlined in~\cite{venkovic2023}.
\end{revenv}

\subsubsection{Algorithm performance in a fixed setting}
\label{subsec:algperf}
Table~\ref{tab:alltestdata} shows various results for the shape problem.
In Table~\ref{tab:midtheta}, we find the results for a medium value for the shape variation parameter $\vartheta$, and we can clearly see that the dimensionality of the problem has a large effect on the algorithm performance.
As expected, the gray-box GPR takes longer to train as the parameter dimension increases.
Additionally, $t_{l-al}$ increases rapidly with the parameter dimension and this increase is more pronounced than in the training time.
This is because the GPR training accounts for parameter dimension anisotropy via the Gaussian process correlation lengths, while the location-allocation does not.
On the other hand, the number of preconditioners placed is relatively stable across parameter dimensions.

Finally, when we compute the total computation times $t_{tot}=t_{train} + t_{l-al} + t_{exec}$, we observe that they range from $8,000$ to $14,000$ seconds, increasing with the parameter dimension.
In these experiments, a single LU preconditioner takes around $t_{pc}=75$ seconds, leading to a baseline computation time of $|W|*t_{pc}=75,000$ seconds when not using a preconditioning strategy, significantly more than the observed range for $t_{tot}$.
Therefore, we have realized a significant reduction of the total computation time.
A single mean based preconditioner does not work very well in this case, as the number of GMRES iterations rises sharply as $\left\| \y \right\|\to\infty$, for which we do not compare against it.
We note that we only know this rises sharply because we have trained the surrogate $m(\cdot)$, it is not clear a-priori whether this is the case or not.
Moreover, we observed values for $m_{\max}$ in Equation~\eqref{eq:f-acq-def} of around $100$ throughout these and other experiments.

    \subsubsection{Varying the magnitude of the shape variations}
    \label{subsec:magnitude-of-shape-variations}
In Table~\ref{tab:hightheta}, the maximal shape variation parameter $\vartheta$ is double from the baseline, just large enough that no self-intersecting curves can occur.
In contrast, Table~\ref{tab:lowtheta} shows the case when the maximal shape variation is halved with respect to the baseline, resulting in significantly smaller shape variations.
This affects matrix values but preserves the overall structure.

Through these experiments, we observe that increasing $\vartheta$ raises training time, while reducing $\vartheta$ lowers it.
For $\vartheta = \frac{1}{4}\vartheta_{\max}(\alpha=2)$, we observe the same behavior in the training time with respect to the parameter dimension: larger parameter dimensions require more training time.
The effect of the parameter dimension is lower in the case of smaller $\vartheta$, as fewer training points are needed to satisfy the stopping criterion.
For $\vartheta = \vartheta_{\max}(\alpha=2)$, we observe similar behavior for $N$, as we did for $\vartheta\frac{1}{2}\vartheta_{\max}(\alpha=2)$ where the training time grows rapidly with the parameter dimension.

The number of preconditioners scales with $\vartheta$ as $m(\y-\yh)$ grows faster with $\|\y-\yh\|$.
Since $t_{pc}$ and $t_{GMRES}$, and therefore, by equation~\eqref{eq:Nratio}, $N_{ratio}$, are independent of $\vartheta$, the `radius of influence' of a preconditioner is reduced, and the preconditioners will need to be packed more densely.

The total computation time is mainly determined by the execution time $t_{exec}$, which increases for larger $\vartheta$.
This is partly due to the different number of preconditioners, and partly due to the slightly larger number of Krylov iterations required due to the larger variations in the resulting matrices.

\subsubsection{Larger anisotropy of the parameter domain}
\label{subsec:anisotropytest}
To study the effects of dimension anisotropy, we modify the decay parameter $\alpha$ in Table~\ref{tab:alpha3}.
Increasing $\alpha$ has two effects: on one hand it decreases the importance of higher parameter dimensions, see equation~\eqref{eq:shapew1infnormdef} and on the other hand, $\vartheta_{\max}(\alpha=3) > \vartheta_{\max}(\alpha=2)$ which we have put to represent total shape variations of the same magnitude.

Foremost, we observe that the values of $t_{train}$ are larger than the values for $\alpha=2$ due to the larger value of $\vartheta$.
On the other side, we see that the training time does not increase much for larger parameter dimensions, clearly showing the effects of the faster decay in the parameter importance.
By contrast, location-allocation cost $t_{l-al}$ still rises sharply since it does not account for anisotropy.
\begin{newenv}
\begin{table}
    \centering
    \begin{tabular}{r||l|l|l|l|l||l|l|l|l|l}
        $N$&\multicolumn{5}{c||}{$\vartheta=\vartheta_{\max}(\alpha=2)$}&\multicolumn{5}{c}{$\vartheta=\frac{1}{2}\vartheta_{\max}(\alpha=2)$}\\\hline
        &$t_{train}$&$t_{l-al}$&$t_{exec}$&$N_{pc}$&$it_{av}$&$t_{train}$&$t_{l-al}$&$t_{exec}$&$N_{pc}$&$it_{av}$\\\hhline{=#=|=|=|=|=#=|=|=|=|=}
        2 & 769 & 9 & 10,355 & 44 & 8 & 402 & 6 & 7,254 & 34 & 6\\
        5 & 2,161 & 145 & 16,577 & 65 & 14 & 515 & 51 & 9,165 & 34 & 9\\
        10 & 4,208 & 1,093 & 18,111 & 69 & 16 & 598 & 167 & 9,287 & 31 & 10\\
        15 & 3,416 & 1,514 & 18,602 & 63 & 16 & 566 & 393 & 9,517 & 31 & 10\\
        20 & 3,605 & 2,994 & 18,803 & 75 & 16 & 635 & 588 & 9,302 & 28 & 11\\
        25 & 5,792 & 5,872 & 19,281 & 67 & 17 & 555 & 943 & 9,693 & 29 & 11\\
        \multicolumn{1}{c}{}\vspace{-10pt}&\multicolumn{5}{c}{\begin{minipage}{4cm}
                                                    \centering
                                                    \subfloat[]{
                                                        \,
                                                        \label{tab:hightheta}
                                                    }
        \end{minipage}}&\multicolumn{5}{c}{\begin{minipage}{5cm}
                                               \centering
                                               \subfloat[][]{
                                                   \label{tab:midtheta}\,\,\,\,
                                               }
        \end{minipage}}\\
        $N$&\multicolumn{5}{c||}{$\vartheta=\frac{1}{4}\vartheta_{\max}(\alpha=2)$}&\multicolumn{5}{c}{$\vartheta=\frac{1}{2}\vartheta_{\max}(\alpha=3)$}\\\hline
        &$t_{train}$&$t_{l-al}$&$t_{exec}$&$N_{pc}$&$it_{av}$&$t_{train}$&$t_{l-al}$&$t_{exec}$&$N_{pc}$&$it_{av}$\\\hhline{=#=|=|=|=|=#=|=|=|=|=}
        2 & 206 & 4 & 5,877 & 19 & 5 & 452 & 6 & 6,996 & 31&6\\
        5 & 204 & 22 & 7,101 & 20 & 7 & 483 & 32 & 8,045 & 29&8\\
        10 & 215 & 85 & 7,124 & 19 & 7 & 471 & 127 & 8,059 & 30&9\\
        15 & 210 & 137 & 7,099 & 17 & 7 & 447 & 218 & 8,232 & 31&8\\
        20 & 233 & 238 & 7,063 & 18 & 7 & 463 & 373 & 7,808 & 29&8\\
        25 & 237 & 408 & 7,203 & 15 & 8 & 497 & 606 & 8,183 & 30&8\\
        \multicolumn{1}{c}{}&\multicolumn{5}{c}{\begin{minipage}{4cm}
                                                    \centering
                                                    \subfloat[]{
                                                        \,
                                                        \label{tab:lowtheta}
                                                    }
        \end{minipage}}&\multicolumn{5}{c}{\begin{minipage}{4cm}
                                               \centering
                                               \subfloat[][]{
                                                   \label{tab:alpha3}\,\,\,\,
                                               }
                                               \clearpage
        \end{minipage}}
    \end{tabular}
    \vspace{-20pt}
    \caption{Results of running Algorithm~\ref{alg:location-allocation} and~\ref{alg:gray-box-gpr} with the shape problem introduced in Section~\ref{subsec:parameterized-domain} with a wavenumber of $k_0=60$ resulting in a mesh with $663,846$ degrees of freedom. Times in seconds and averaged over multiple runs.}
    \label{tab:alltestdata}
\end{table}
\end{newenv}
\subsubsection{Lower and higher frequency}
\label{subsec:freqtest}
To change problem conditioning and matrix sizes, we modify the frequency $k_0$ of the Helmholtz equation, shown in Table~\ref{tab:shapefrequencies}.
We consider two cases: a lower frequency of $k_0=30$ and a higher frequency $k_0=120$.
To compare the numerical performance, we keep the size of $W$ the same at 1000 points.
However, we emphasise that, when one is interested in the error in a surrogate or an average using Monte Carlo or sparse grid methods, the number of points should be adjusted with $k_0$~\cite{hiptmair2024}.
Due to computational limits, we run our algorithm once for $k_0=120$, without averaging over multiple runs.

We observe a large difference in the overall magnitude of the computation times: a low value of $k_0$ corresponds to low computation times and vice versa in the high frequency case.
In the case $k_0=30$, we compare the total computation times with the baseline $|W|*t_{pc}\approx 1000*1.5=1500$ and conclude that we have improved the computational load, but barely when the parameter dimension is large.
In the case $k_0=120$, we have $|W|*t_{pc}\approx 1000*1000s=1,000,000s\approx11.5\text{ days}$, and we have improved the computational load by a factor of 10.

In the case $k_0=30$, the location-allocation cost does not increase as much as the case $k_0=60$ (Table~\ref{tab:midtheta}), although the number of preconditioners is roughly the same as the dimension increases.
This is in contrast to the case $k_0=120$, where $t_{l-al}$ increases significantly.
We have seen similar increases before in Table~\ref{tab:alltestdata}, where the increase was caused by the higher number of placed preconditioners.
In this case, the values of $N_{pc}$ are similar and the increase is due to additional location-allocation steps.
The algorithm runs until no further improvements are found or iteration costs exceed the improvement gained, and higher computation times allow more resources for location-allocation iterations.
The lower values of $t_{l-al}$ in the case $k_0=30$ can be explained similarly.
Finally, the number of preconditioners is higher for $k_0=30$, which is because $\tau_{pc}$ decreases faster than $\tau_{Krylov}$ as the matrix dimension shrinks.

\begin{table}
    \centering
    \begin{tabular}{r||l|l|l|l|l||l|l|l|l|l}
        $N$&\multicolumn{5}{c||}{$k_0=30$, $\alpha=2$}&\multicolumn{5}{c}{$k_0=120$, $\alpha=2$}\\\hline
        &$t_{train}$&$t_{l-al}$&$t_{exec}$&$N_{pc}$&$it_{av}$&$t_{train}$&$t_{l-al}$&$t_{exec}$&$N_{pc}$&$it_{av}$\\\hhline{=#=|=|=|=|=#=|=|=|=|=}
        2 & 14 & 7 & 243 & 49 & 4 & 3,291 & 6 & 46,364 & 23 & 9\\
        5 & 14 & 40 & 339 & 60 & 6 & 10,848 & 60 & 71,575 & 26 & 17\\
        10 & 20 & 125 & 362 & 57 & 7 & 22,079 & 458 & 81,134 & 24 & 20\\
        15 & 20 & 224 & 359 & 48 & 7 & 29,513 & 1,169 & 77,540 & 22 & 19\\
        20 & 36 & 254 & 346 & 30 & 7 & 21,833 & 1,676 & 77,496 & 18 & 22\\
        25 & 31 & 349 & 349 & 26 & 8 & 31,626 & 3,010 & 77,043 & 22 & 21\\
        \multicolumn{1}{c}{}&\multicolumn{5}{c}{\begin{minipage}{4cm}
                                                    \centering
                                                    \subfloat[]{
                                                        \,
                                                        \label{tab:k30}
                                                    }
        \end{minipage}}&\multicolumn{5}{c}{\begin{minipage}{4cm}
                                                    \centering
                                                    \subfloat[][]{
                                                        \label{tab:k120}\,\,\,\,
                                                    }
        \end{minipage}}
    \end{tabular}
    \vspace{-20pt}
    \caption{Test results for the shape expansion with $\alpha=2$, $\vartheta=2^{-1}\vartheta_{\max}(\alpha=2)$ and frequencies $k_0=30$ and $k_0=120$ resulting in meshes of $84,255$ and $5,276,457$ degrees of freedom, respectively. Times in seconds and averaged over multiple runs.}
    \label{tab:shapefrequencies}
\end{table}

\subsection{Algorithm performance for the affine expansion}\label{subsec:algorithm-performance-for-the-affine-expansion}
Table~\ref{tab:eta025} shows results for the case of affine expansion outlined in Section~\ref{subsec:disjoint-density}, with a fully isotropic parameter domain.
This case helps us to study the effects of the parameter anisotropy.
We present results with $\eta_i=\frac{1}{4}$, for $i=1\cdots N$.
This value is not very high, which is reflected in the results.

The most obvious effect is shown in the column $N_{pc}$, showing only one preconditioner for higher parameter dimensions.
This is due to the concentration effect (Remark~\ref{rem:high-dimensional-parameter-space}) and is in stark contrast to the anisotropic cases discussed before.
Thus, we are in fact using mean-based preconditioning, with additional training overhead $t_{train}$.

Additionally, the training time $t_{train}$ is not significantly impacted by the parameter dimension.
This is because the prior mean approximates the surrogate very well, requiring only a few training points.
To examine this further, we approximate the root-mean-square error (RMSE) for $N=2,10,25$ during training.
Fig.~\ref{fig:pietrain} shows that the GPR trains fairly well in the lower dimensional limit, decreasing the RMSE from 10 to 6 quickly, after which the training stagnates.
For high dimensions ($N=25$), the prior fits the model very well, but training does not improve the surrogate.
This highlights the curse of dimensionality, showing that parameter anisotropy is necessary for optimal performance.
\begin{figure}
    \begin{minipage}[t]{.50\linewidth}
        \vspace{10pt}
        \begin{tabular}[t]{r|l|l|l|l|l}
            $N$&$t_{train}$&$t_{l-al}$&$t_{exec}$&$N_{pc}$&$it_{av}$\\\hline
            2 & 1,551 & 3 & 5,339 & 12 & 5 \\
            5 & 1,517 & 11 & 7,329 & 7 & 9 \\
            10 & 1,756 & 0 & 7,167 & 1 & 10 \\
            15 & 1,627 & 0 & 7,501 & 1 & 9 \\
            20 & 1,616 & 0 & 7,383 & 1 & 9 \\
            25 & 1,756 & 0 & 7,168 & 1 & 9 \\
        \end{tabular}
        \captionof{table}{Tabulated results in seconds for affine expansion with $\eta_i=\frac{1}{4}$.}
        \label{tab:eta025}
        \end{minipage}\hfill%
    \begin{minipage}[t]{.48\linewidth}
        \vspace{0pt}
        \begin{tikzpicture}
            \begin{axis}[
                ymin=0,
                ymax=10,
                xmin=0,
                xmax=200,
                width=\textwidth,
                height=0.725\textwidth,
                ylabel=RMSE,
                ylabel style={font=\small},
                ylabel near ticks
            ]
                \addplot[black,thick] table [mark=none, x=iter, y=d2]{data/pietrain025.dat};
                \addlegendentry{$N\eq 2$}
                \addplot[black,thick,dotted] table [mark=none, x=iter, y=d10]{data/pietrain025.dat};
                \addlegendentry{$N\eq 10$}
                \addplot[black,thick,dashed] table [mark=none, x=iter, y=d25]{data/pietrain025.dat};
                \addlegendentry{$N\eq 25$}
            \end{axis}
        \end{tikzpicture}
        \caption{Training RMSE for affine expansion.}
        \label{fig:pietrain}
    \end{minipage}
\end{figure}

\subsection{Preconditioner placement}\label{subsec:preconditioner-placement}
\begin{figure}
    \centering
    \begin{subfigure}{0.50\textwidth}
        \includegraphics[width=\linewidth]{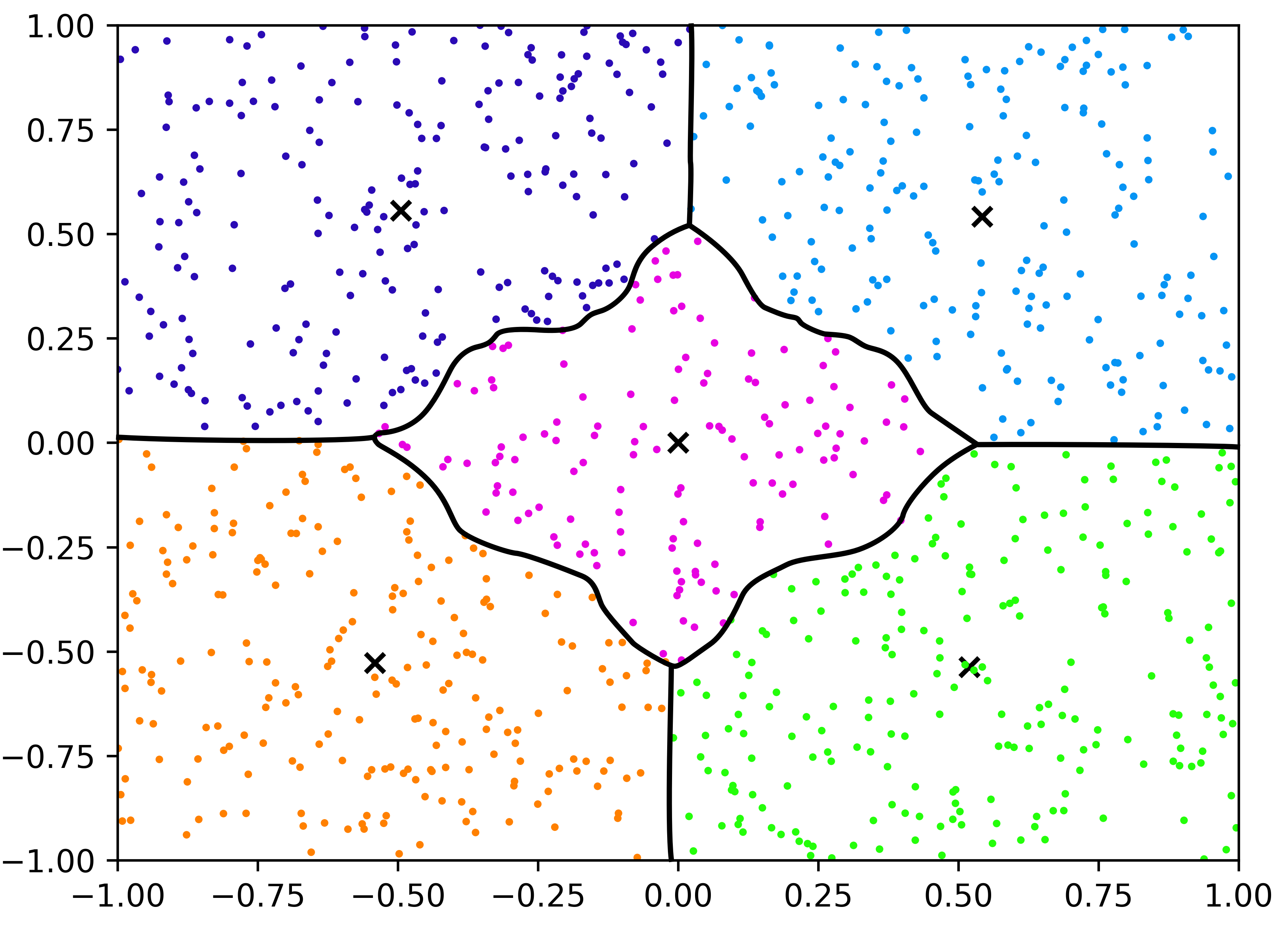}
        \caption{} \label{fig:1a}
    \end{subfigure}%
    \begin{subfigure}{0.50\textwidth}
        \includegraphics[width=\linewidth]{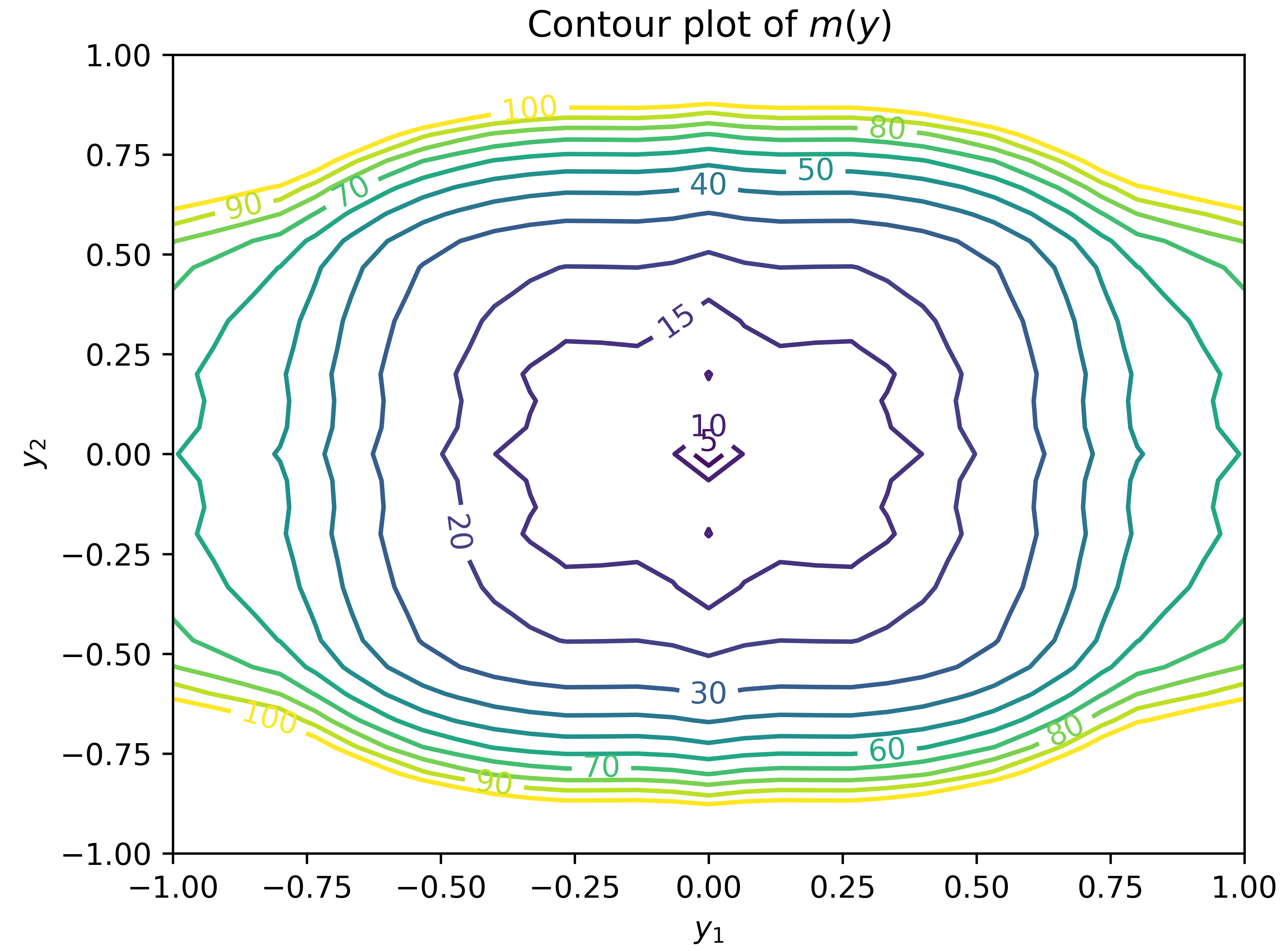}
        \caption{} \label{fig:1b}
    \end{subfigure}%
    \vspace*{-15pt}
    \caption{Result of the location-allocation algorithm; Fig.~\ref{fig:1a} shows the preconditioner attributions together with the underlying partition of $W$, and Fig.~\ref{fig:1b} shows a contour plot of the underlying distance function $m(\cdot)$.}
    \label{fig:voronoi_diagram}
\end{figure}

We analyze the efficacy of Algorithm~\ref{alg:location-allocation} in placing the preconditioners in Fig.~\ref{fig:voronoi_diagram}.
To make this insightful, we use a low parameter dimension of $N=2$, which allows for clear visualizations, and we sample uniformly over the parameter space.
Fig.~\ref{fig:1a} shows the result of the location-allocation algorithm.
The colored dots represent the parameter values, and they are grouped by their assigned preconditioner, which is indicated by color.
Moreover, in black, the borders of the generalized Voronoi diagram are shown.

From Fig.~\ref{fig:voronoi_diagram} it is clear that the edges of the Voronoi diagram are curved because the Voronoi distance is based on $m(\cdot)$ shown in Fig.~\ref{fig:1b}, which deviates from the actual distance to the preconditioner.
Moreover, Fig.~\ref{fig:1a} shows that not all preconditioners are placed at parameter locations, see Remark~\ref{rem:pcnotatcol}.

\subsection{Number of preconditioners}\label{subsec:Npc_experiment}
Fig.~\ref{fig:Npc_tikz} shows the efficacy of the preconditioner count selection in Algorithm~\ref{alg:location-allocation}.
As discussed in Section~\ref{subsec:determining-Npc}, we re-use the initialization approach for efficiency.
To compare against the optimal choice, we preform greedy initialization and execute Algorithm~\ref{alg:location-allocation} after each step.
This allows us to check whether our approach is close to the optimal $N_{pc}$.
We perform this comparison for the parameterized shape problem with parameter dimension $N=15$.

The red line represents the greedy initialization costs, and the blue line corresponds to the costs after the location-allocation are steps taken as long as it is worth the computation time.
As expected, there is an optimal number of preconditioners, as the curves exhibit a minimum.
The red line being above the blue line indicates that the location-allocation procedure improves preconditioner placement.
The black curve shows results if we run location-allocation until full convergence.
Since it overlaps almost perfectly with the blue curve, cutting location-allocation short has little impact, and our strategy is essentially optimal.

\begin{revenv}
    The cost curve in Fig.~\ref{fig:Npc_tikz} decreases fast for low numbers of preconditioners, indicating that using only a few extra preconditioners can impact the total computation cost.
    Therefore, if, due to system requirements, there is a strict limit on computational resources, $N_{pc}$ can be chosen to be less than $N_{pc}^\ast$ and still give performance improvements.
\end{revenv}

\begin{revenv}
    Moreover, we observe that the black and blue curves are not completely smooth.
    This occurs because the location-allocation heuristic ends up in local minima.
    However, the volatility is small, and does not majorly impact the resulting strategy.
\end{revenv}
Finally, we see the three curves converging for large $N_{pc}$, indicating that the greedy initialization performs well for small $\frac{|W|}{N_{pc}}$.

For this run, the chosen value for $N_{pc}$ is at $N_{pc}^*=35$, marked by a vertical line in Fig.~\ref{fig:Npc_tikz}.
This is lower than the minimum from the greedy initialization, which happened around $N_{pc}\approx 45$.
Still, the final strategy remains close to optimal.
Variance from local minima overshadows any improvements, making this a cost-effective method.

\begin{figure}
    \centering
    \npcgraph{data/Npc3mar/Npc_data_10.dat}
    \vspace{-20pt}
    \caption{
Effect of different $N_{pc}$ values for the parameterized shape problem.
The red curve shows the greedy initialization, the blue curve the location-allocation result from Algorithm~\ref{alg:location-allocation}, and the black dashed curve the result of the location-allocation algorithm until convergence.
$N_{pc}^*$ is the value where Algorithm~\ref{alg:location-allocation} terminates.
}
    \label{fig:Npc_tikz}
\end{figure}

\subsection{Convergence of the surrogate model}\label{subsec:gpr_reliability_experiment}
The reliability of our preconditioning strategy hinges on a good estimation $m(\cdot)$ for the number of GMRES iterations, which we will explore in this section by using the parameter domain problem with parameter dimension $N=15$.
The results are shown in Fig.~\ref{fig:surrogate_convergence}.

We expect $m(\cdot)$ to predict iterations more accurately with more training data.
To verify this, we train the GPR with different numbers of training data $N_{train}$.
After each training point we add, we compute the root-mean-square error of $m(\cdot)$ over the domain $\{\y \in Y | m(\y) < N_{ratio} / 2\}$ using a Monte Carlo estimate.
This domain is taken such that we measure the accuracy of the surrogate over the relevant domain, as a parameter value with $m(\y)>\frac{N_{ratio}}{2}$ will most likely get assigned to another preconditioner even further from the origin.
We use a parameter dimension $N=15$, anisotropy $\alpha=2$, and $\vartheta=\frac{1}{2}\vartheta_{\max}(\alpha=2)$.

Fig.~\ref{fig:surrogate_convergence} shows the results: the surrogate accuracy improves with more training data but plateaus at a saturation point.
Next to this, we observe that the disagree ratio fluctuates heavily, emphasizing the need for smoothing.
Once smoothened, the disagree ratio decays as we add more training data, and it crosses the $1\%$ threshold before the saturation threshold.
We observe this behavior consistently across multiple runs and varying parameter dimensions, although the exact stopping point varies.
Despite the aggressiveness of the stopping criterion, the surrogate still enables effective preconditioner placement.

Finally, the RMSE converges to a value of approximately 3.5.
This value reflects the modeling error between the actual number of GMRES iterations needed and our surrogate $m(\cdot)$.
This originates from our kernel choice, where we assumed a dimensional splitting, together with the symmetry around the origin.
These assumptions reduce accuracy but mitigate the curse of dimensionality.
A RMSE of 3.5 is still very acceptable, and the gained computational improvements are well worth it.
\begin{figure}
    \centering
    \begin{tikzpicture}
        \begin{axis}[
            label style={font=\small},
            ticklabel style = {font=\small},
            axis y line=left,
            ymax=12,
            ymin=3,
            xmin=0,
            xmax=100,
            ytick={10, 8, 6, 4, 2},
            yticklabels={10, 8, 6, 4, 2},
            width=0.90\textwidth,
            height=0.35\textwidth,
            ylabel near ticks,
            legend columns=4,
            ylabel=RMSE]
            \addplot[black,thick] table [mark=none, x=Ntrain, y=Dimsum15]{data/LOOCV/SP_retest.dat};
            \label{plot_one}
        \end{axis}
        \begin{axis}[
            axis y line=right,
            axis x line=none,
            ytick={1, 0.1, 0.01, 0.001, 0.0001},
            ymode=log,
            legend style={font=\fontsize{8}{7}\selectfont},
            label style={font=\small},
            ticklabel style = {font=\small},
            ymax=12,
            xmin=0,
            xmax=100,
            width=0.90\textwidth,
            height=0.35\textwidth,
            legend columns=4,
            ylabel=disagree ratio]

            \addlegendimage{/pgfplots/refstyle=plot_one}\addlegendentry{RMSE}
            \addplot[black,thick,dashed] table [mark=none, x=Ntrain, y=SP15]{data/LOOCV/SP_retest.dat};
            \addlegendentry{disagree ratio}
            \addplot[black,thick,densely dotted, line width=1pt, color=red] table [mark=none, x=Ntrain, y=Smooth15]{data/LOOCV/SP_retest.dat};
            \addlegendentry{smooth disagree ratio}
            \addplot[color=gray, line width=1pt] coordinates {(0, 0.01) (200, 0.01)};
            \addlegendentry{$1\%$}
        \end{axis}
    \end{tikzpicture}
    \caption{Comparison of the RMSE and disagree ratio for the stabilizing predictions (SP) stopping criterion.
    The parameter dimension is $N=15$ with $\vartheta=\frac{1}{2}\vartheta_{\max}$ and $\alpha=2$.
    The solid line represents the RMSE (left axis), the dashed line the disagree ratio, and the dotted line a trailing average of the disagree ratio (right axis).
    The solid gray line marks the $1\%$ threshold in the stopping criterion.
    }\label{fig:surrogate_convergence}
\end{figure}

    \section{Discussion and conclusion}\label{sec:conclusions}
    In this work, we developed a strategy to place multiple preconditioners in the parameter space that can reduce the computational cost of repeatedly solving parameterized linear systems by an order of magnitude.
We use a two-step process: first, we learn a surrogate for the number of Krylov iterations with gray-box Gaussian process regression trained with a cost-aware active learning strategy.
This surrogate model allows for estimating optimal number of preconditioners using a greedy approach and finally use a location-allocation algorithm to optimize their locations.
We have applied and studied the algorithm thoroughly using a Helmholtz scattering problem.

To choose the \textit{Gaussian process prior}, we use a-priori bounds on the number of GMRES iterations, which are available for the parametric Helmholtz equation~\cite{graham2021}.
The modeling error caused by the upper bounds is corrected by training the GPR and its hyperparameters.
In different applications where such upper bounds are not available, the distance to the origin could be used for the prior mean.
Because we use a-priori bounds this way, the prior is more informative and less training is required.

For the \textit{Gaussian process kernel}, we have used a symmetrized Matérn kernel per dimension, summed over all dimensions.
This assumption limits the richness of the function space we approximate, and numerical experiments have shown that using a richer full Matérn kernel allows the surrogate to become more accurate, but it suffers more from the curse of dimensionality.
Moreover, we notice that, to find large improvements in the total computational burden, a `good' rather than the most accurate surrogate is enough.

For \textit{placing the preconditioners}, we solve a high-dimensional optimization problem where the dimension scales with the number of preconditioners times the parameter dimension.
We use a location-allocation approach to handle this high dimensionality during preconditioner placement.
Since the iterations are costly, we perform them until the improvement is smaller than the computation time.
The main cost lies in the location step, where the geometric median of each partition cell $W_k$ has to be computed.
Although the cost of this preprocessing step using an off-the-shelf optimizer (L-BFGS-B) is within our requirements, further savings could be achieved with an ad-hoc optimizer.
Unlike previous work~\cite{venkovic2023}, we do not require the partitions $W_k$ to be of equal cardinality, possibly resulting in preconditioners with few allocated parameter locations and others with many.
However, this is not problematic as we have not considered any multithreaded approaches, which might alter this tradeoff.

Determining the optimal \textit{number of preconditioners} is difficult, as we cannot execute the location-allocation algorithm multiple times, and we thus use a greedy approach.
Given the existence of local minima in the preconditioner placement, any further improvements will be overshadowed by these local minima.
Since the number and the placement of preconditioners depend on the ratio of preconditioner computation time to Krylov iteration time, our algorithm is the most reliable when there are either no other computational loads or a constant load.
\begin{revenv}
    Finally, after the training phase has concluded, the algorithm might conclude that a single preconditioner, i.e., mean-based precondioning, is the best strategy.
    In this case, the training phase was still useful, as the training would have already solved quite some linear systems.
    The overhead of the Gaussian Process training is modest in comparison.
\end{revenv}

Because of \textit{concentration of measure} effects for high $N$, our approach will not find significant improvements for high dimensional isotropic parameter spaces.
In this case, the optimal number of preconditioners is either one (mean-based preconditioning) or to place a preconditioner at each parameter location.
Adaptive methods can still determine which case is applicable, ensuring that the correct case is identified.

\begin{revenv}
    The \textit{memory requirements} of our approach are modest because only one preconditioner is needed in system memory at once.
    The mean-based preconditioner is needed in the training phase, and the resulting strategy can be executed sequentially by keeping only one preconditioner in memory at the time.
    When considering, for example, problems with a higher spatial (e.g. $d=3$) dimension, the memory requirements can be accommodated by switching to sparser preconditioning strategies, such as in~\cite{bollhofer2009,engquist2011b,poulson2013}.
    However, the combined efficacy of these methods is still an open problem.

    In this work, we have considered a \emph{LU preconditioner}.
    A similar approach could be explored for different preconditioners.
    However, as the algorithm re-uses the preconditioner multiple times, an expensive preconditioner such as LU will result in more savings, as the price paid to compute it can be divided over multiple solves while the benefits are applicable in each solve separately.
    Moreover, to get the estimate for $m$, we exploited the bounds in~\cite{graham2021} for LU preconditioners.
    To extend the method to different preconditioners, we either need similar estimates for these preconditioners or a more general Gaussian Process without the gray-box approach we have used.
\end{revenv}

Further directions that could be of interest are, among others, \begin{revenv}
                                                                    the employment of GMRES restart,
\end{revenv}the application of different surrogate models, multithreaded approaches, and the use of preconditioners different from the full LU decomposition.

    \bibliographystyle{plain} 
    \bibliography{library}

\begin{thebibliography}{10}

\bibitem{addy2025}
{\sc E.~J. Addy, J.~Latz, and A.~L. Teckentrup}, {\em Lengthscale-informed
  sparse grids for kernel methods in high dimensions}, ArXiv, abs/2506.07797
  (2025), \url{https://api.semanticscholar.org/CorpusID:279250459}.

\bibitem{aggarwal1973}
{\sc C.~C. Aggarwal, A.~Hinneburg, and D.~A. Keim}, {\em On the {{Surprising
  Behavior}} of {{Distance Metrics}} in {{High Dimensional Space}}}, in
  Database {{Theory}} -- {{ICDT}} 2001, Lecture {{Notes}} in {{Computer
  Science}}, Springer, 1973, pp.~420--435.

\bibitem{alnaes2014a}
{\sc M.~S. Aln{\ae}s, A.~Logg, K.~B. {\O}lgaard, M.~E. Rognes, and G.~N.
  Wells}, {\em Unified form language: {{A}} domain-specific language for weak
  formulations of partial differential equations}, ACM Transactions on
  Mathematical Software, 40 (2014), pp.~1--37,
  \url{https://doi.org/10.1145/2566630}.

\bibitem{astudillo2022}
{\sc R.~Astudillo and P.~I. Frazier}, {\em Thinking inside the box: {{A}}
  tutorial on grey-box {{Bayesian}} optimization}, in {{WSC}} '21:
  {{Proceedings}} of the {{Winter Simulation Conference}}, Jan. 2022,
  \url{https://doi.org/10.1109/WSC52266.2021.9715343}.

\bibitem{azulay2022multigrid}
{\sc Y.~Azulay and E.~Treister}, {\em Multigrid-augmented deep learning
  preconditioners for the {{Helmholtz}} equation}, SIAM Journal on Scientific
  Computing, 45 (2022), pp.~S127--S151.

\bibitem{babuska1997}
{\sc I.~M. Babu{\v s}ka and S.~A. Sauter}, {\em Is the {{Pollution Effect}} of
  the {{FEM Avoidable}} for the {{Helmholtz Equation Considering High Wave
  Numbers}}?}, SIAM Journal on Numerical Analysis, 34 (1997), pp.~2392--2423,
  \url{https://doi.org/10.1137/S0036142994269186}.

\bibitem{baratta2023}
{\sc I.~A. Baratta, J.~P. Dean, J.~S. Dokken, M.~Habera, J.~S. Hale, C.~N.
  Richardson, M.~E. Rognes, M.~W. Scroggs, N.~Sime, and G.~N. Wells}, {\em
  {{DOLFINx}}: {{The}} next generation {{FEniCS}} problem solving environment},
  Dec. 2023, \url{https://doi.org/10.5281/ZENODO.10447666}.

\bibitem{binois2022}
{\sc M.~Binois and N.~Wycoff}, {\em A {{Survey}} on {{High-dimensional Gaussian
  Process Modeling}} with {{Application}} to {{Bayesian Optimization}}}, ACM
  Transactions on Evolutionary Learning and Optimization, 2 (2022), pp.~1--26,
  \url{https://doi.org/10.1145/3545611}.

\bibitem{bishop2006}
{\sc C.~M. Bishop}, {\em Pattern Recognition and Machine Learning}, Information
  Science and Statistics, Springer, New York, 2006.

\bibitem{bloodgood2009}
{\sc M.~Bloodgood and V.~K. Shanker}, {\em A {{Method}} for {{Stopping Active
  Learning Based}} on {{Stabilizing Predictions}} and the {{Need}} for
  {{User-Adjustable Stopping}}}, in Proceedings of the {{Thirteenth
  Conference}} on {{Computational Natural Language Learning}} ({{CoNLL-2009}}),
  Boulder, Colorado, 2009, Association for Computational Linguistics,
  pp.~39--47.

\bibitem{bollhofer2009}
{\sc M.~Bollh{\"o}fer, M.~J. Grote, and O.~Schenk}, {\em Algebraic {{Multilevel
  Preconditioner}} for the {{Helmholtz Equation}} in {{Heterogeneous Media}}},
  SIAM Journal on Scientific Computing, 31 (2009), pp.~3781--3805,
  \url{https://doi.org/10.1137/080725702}.

\bibitem{bora2025}
{\sc D.~J. Bora and M.~K. Mishra}, {\em Comparative {{Evaluation}} of {{Hard}}
  and {{Soft Clustering}} for {{Precise Brain Tumor Segmentation}} in {{MR
  Imaging}}}, Journal of Advances in Mathematics and Computer Science, 40
  (2025), pp.~127--141, \url{https://doi.org/10.9734/jamcs/2025/v40i92050},
  \url{https://arxiv.org/abs/2509.05340}.

\bibitem{brimberg2008}
{\sc J.~Brimberg, P.~Hansen, N.~Mladonovic, and S.~Salhi}, {\em A survey of
  solution methods for the continuous location allocation problem},
  International Journal of Operational Research, 5 (2008), pp.~1--12.

\bibitem{brown2022}
{\sc S.~B. Brown, S.~Abhyankar, M.~F. Adams, S.~Benson, J.~Dener, P.~Brune,
  K.~Buschelman, E.~M. Constantinescu, L.~Dalcin, and A.~Jolivet}, {\em {{PETSc
  Web}} page}, 2022.

\bibitem{castrillon-candas2016}
{\sc J.~E. {Castrill{\'o}n-Cand{\'a}s}, F.~Nobile, and R.~F. Tempone}, {\em
  Analytic regularity and collocation approximation for elliptic {{PDEs}} with
  random domain deformations}, Computers and Mathematics with Applications, 71
  (2016), pp.~1173--1197, \url{https://doi.org/10.1016/j.camwa.2016.01.005}.

\bibitem{chandler-wilde2008}
{\sc S.~N. {Chandler-Wilde} and P.~Monk}, {\em Wave-{{Number-Explicit Bounds}}
  in {{Time-Harmonic Scattering}}}, SIAM Journal on Mathematical Analysis, 39
  (2008), pp.~1428--1455, \url{https://doi.org/10.1137/060662575}.

\bibitem{chaumont-frelet2023}
{\sc T.~{Chaumont-Frelet}, A.~Moiola, and E.~A. Spence}, {\em Explicit bounds
  for the high-frequency time-harmonic {{Maxwell}} equations in heterogeneous
  media}, Journal de Math\'ematiques Pures et Appliqu\'ees, 179 (2023),
  pp.~183--218, \url{https://doi.org/10.1016/j.matpur.2023.09.004}.

\bibitem{chen2022}
{\sc Y.~Chen, B.~Dong, and J.~Xu}, {\em Meta-{{MgNet}}: {{Meta}} multigrid
  networks for solving parameterized partial differential equations}, Journal
  of Computational Physics, 455 (2022), p.~110996,
  \url{https://doi.org/10.1016/j.jcp.2022.110996}.

\bibitem{chew1985}
{\sc P.~Chew and R.~L. Dyrsdale}, {\em Voronoi diagrams based on convex
  distance functions}, in Proceedings of the First Annual Symposium on
  {{Computational}} Geometry - {{SCG}} '85, Baltimore, Maryland, United States,
  1985, ACM Press, pp.~235--244, \url{https://doi.org/10.1145/323233.323264}.

\bibitem{church2022}
{\sc R.~L. Church, Z.~Drezner, and A.~Tamir}, {\em Extensions to the {{Weber}}
  problem}, Computers \& Operations Research, 143 (2022), p.~105786,
  \url{https://doi.org/10.1016/j.cor.2022.105786}.

\bibitem{cohen2018}
{\sc A.~Cohen, C.~Schwab, and J.~Zech}, {\em Shape {{Holomorphy}} of the
  {{Stationary Navier--Stokes Equations}}}, SIAM Journal on Mathematical
  Analysis, 50 (2018), pp.~1720--1752,
  \url{https://doi.org/10.1137/16M1099406}.

\bibitem{cohen1960}
{\sc J.~Cohen}, {\em A {{Coefficient}} of {{Agreement}} for {{Nominal
  Scales}}}, Educational and Psychological Measurement, 20 (1960), pp.~37--46,
  \url{https://doi.org/10.1177/001316446002000104}.

\bibitem{contreras2018}
{\sc A.~A. Contreras, P.~Mycek, O.~P. Le~Ma{\^i}tre, F.~Rizzi, B.~Debusschere,
  and O.~M. Knio}, {\em Parallel {{Domain Decomposition Strategies}} for
  {{Stochastic Elliptic Equations Part B}}: {{Accelerated Monte Carlo
  Sampling}} with {{Local PC Expansions}}}, SIAM Journal on Scientific
  Computing, 40 (2018), pp.~C547--C580,
  \url{https://doi.org/10.1137/17M1132197}.

\bibitem{cooper1963}
{\sc L.~Cooper}, {\em Location-{{Allocation Problems}}}, Operations Research,
  11 (1963), pp.~331--343, \url{https://doi.org/10.1287/opre.11.3.331}.

\bibitem{cooper1981}
{\sc L.~Cooper and N.~Katz}, {\em The {{Weber}} problem revisited}, Computers
  \& Mathematics with Applications, 7 (1981), pp.~225--234,
  \url{https://doi.org/10.1016/0898-1221(81)90082-1}.

\bibitem{cui2025}
{\sc C.~Cui, K.~Jiang, and S.~Shu}, {\em A {{Neural Multigrid Solver}} for
  {{Helmholtz Equations}} with {{High Wavenumber}} and {{Heterogeneous
  Media}}}, SIAM Journal on Scientific Computing, 47 (2025), pp.~C655--C679,
  \url{https://doi.org/10.1137/24M1654397}.

\bibitem{drezner2009}
{\sc Z.~Drezner}, {\em On the convergence of the generalized {{Weiszfeld}}
  algorithm}, Annals of Operations Research, 167 (2009), pp.~327--336,
  \url{https://doi.org/10.1007/s10479-008-0336-z}.

\bibitem{duvenaud2014}
{\sc D.~Duvenaud}, {\em Automatic Model Construction with {{Gaussian}}
  Processes}, PhD thesis, Apollo - University of Cambridge Repository, Nov.
  2014, \url{https://doi.org/10.17863/CAM.14087}.

\bibitem{duvenaud2011}
{\sc D.~Duvenaud, H.~Nickisch, and C.~E. Rasmussen}, {\em Additive {{Gaussian
  Processes}}}, in Advances in {{Neural Information Processing Systems}},
  vol.~24, Curran Associates, Inc., 2011.

\bibitem{eckhardt1980}
{\sc U.~Eckhardt}, {\em Weber's problem and {{Weiszfeld}}'s algorithm in
  general spaces}, Mathematical Programming, 18 (1980), pp.~186--196,
  \url{https://doi.org/10.1007/BF01588313}.

\bibitem{eiermann2007}
{\sc M.~Eiermann, O.~G. Ernst, and E.~Ullmann}, {\em Computational aspects of
  the stochastic finite element method}, Computing and Visualization in
  Science, 10 (2007), pp.~3--15,
  \url{https://doi.org/10.1007/s00791-006-0047-4}.

\bibitem{elman1982}
{\sc H.~C. Elman}, {\em Iterative {{Methods}} for {{Large}}, {{Sparse}},
  {{Nonsymmetric Systems}} of {{Linear Equations}}}, PhD thesis, Yale
  University, Apr. 1982.

\bibitem{embree1999}
{\sc M.~Embree}, {\em How descriptive are {{GMRES}} convergence bounds?}, NA
  Report 99, 8 (1999).

\bibitem{engquist2011b}
{\sc B.~Engquist and L.~Ying}, {\em Sweeping {{Preconditioner}} for the
  {{Helmholtz Equation}}: {{Moving Perfectly Matched Layers}}}, Multiscale
  Modeling \& Simulation, 9 (2011), pp.~686--710,
  \url{https://doi.org/10.1137/100804644}.

\bibitem{ernst2009}
{\sc O.~G. Ernst, C.~E. Powell, D.~J. Silvester, and E.~Ullmann}, {\em
  Efficient {{Solvers}} for a {{Linear Stochastic Galerkin Mixed Formulation}}
  of {{Diffusion Problems}} with {{Random Data}}}, SIAM Journal on Scientific
  Computing, 31 (2009), pp.~1424--1447,
  \url{https://doi.org/10.1137/070705817}.

\bibitem{geuzaine2009}
{\sc C.~Geuzaine and J.-F. Remacle}, {\em Gmsh: {{A}} 3-{{D}} finite element
  mesh generator with built-in pre- and post-processing facilities},
  International Journal for Numerical Methods in Engineering, 79 (2009),
  pp.~1309--1331, \url{https://doi.org/10.1002/nme.2579}.

\bibitem{ghanem1996}
{\sc R.~G. Ghanem and R.~M. Kruger}, {\em Numerical solution of spectral
  stochastic finite element systems}, Computer Methods in Applied Mechanics and
  Engineering, 129 (1996), pp.~289--303,
  \url{https://doi.org/10.1016/0045-7825(95)00909-4}.

\bibitem{ginsbourger2010}
{\sc D.~Ginsbourger, R.~Le~Riche, and L.~Carraro}, {\em Kriging {{Is
  Well-Suited}} to {{Parallelize Optimization}}}, in Computational
  {{Intelligence}} in {{Expensive Optimization Problems}}, L.~M. Hiot, Y.~S.
  Ong, Y.~Tenne, and C.-K. Goh, eds., vol.~2, Springer Berlin Heidelberg,
  Berlin, Heidelberg, 2010, pp.~131--162,
  \url{https://doi.org/10.1007/978-3-642-10701-6_6}.

\bibitem{giraud:hal-05157038}
{\sc L.~Giraud, C.~Kruse, P.~Mycek, M.~Shpakovych, and Y.~Xiang}, {\em Neural
  network preconditioning: A case study for the solution of the parametric
  {{Helmholtz}} equation}, Tech. Report RR-9593, Inria Centre at the University
  of Bordeaux, France, July 2025.

\bibitem{graham2019}
{\sc I.~G. Graham, O.~R. Pembery, and E.~A. Spence}, {\em The {{Helmholtz}}
  equation in heterogeneous media: {{A}} priori bounds, well-posedness, and
  resonances}, Journal of Differential Equations, 266 (2019), pp.~2869--2923,
  \url{https://doi.org/10.1016/j.jde.2018.08.048}.

\bibitem{graham2021}
{\sc I.~G. Graham, O.~R. Pembery, and E.~A. Spence}, {\em Analysis of a
  {{Helmholtz}} preconditioning problem motivated by uncertainty
  quantification}, Advances in Computational Mathematics, 47 (2021), pp.~1--39,
  \url{https://doi.org/10.1007/s10444-021-09889-0}.

\bibitem{harbrecht2016}
{\sc H.~Harbrecht, M.~Peters, and M.~Siebenmorgen}, {\em Analysis of the domain
  mapping method for elliptic diffusion problems on random domains}, Numerische
  Mathematik, 134 (2016), pp.~823--856,
  \url{https://doi.org/10.1007/s00211-016-0791-4}.

\bibitem{hiptmair2018}
{\sc R.~Hiptmair, L.~Scarabosio, C.~Schillings, and C.~Schwab}, {\em Large
  deformation shape uncertainty quantification in acoustic scattering},
  Advances in Computational Mathematics, 44 (2018), pp.~1475--1518,
  \url{https://doi.org/10.1007/s10444-018-9594-8}.

\bibitem{hiptmair2024}
{\sc R.~Hiptmair, C.~Schwab, and E.~A. Spence}, {\em Frequency-explicit shape
  holomorphy in uncertainty quantification for acoustic scattering}, SIAM/ASA
  Journal on Uncertainty Quantification, 13 (2025), pp.~1904--1949,
  \url{https://doi.org/10.1137/24M1688643},
  \url{https://doi.org/10.1137/24M1688643},
  \url{https://arxiv.org/abs/https://doi.org/10.1137/24M1688643}.

\bibitem{horn1991}
{\sc R.~A. Horn and C.~R. Johnson}, {\em Topics in {{Matrix Analysis}}},
  Cambridge University Press, Apr. 1991,
  \url{https://doi.org/10.1017/cbo9780511840371}.

\bibitem{howell2009}
{\sc C.~Howell}, {\em Lifting the curse of dimensionality: {{A}} europeanist's
  perspective}, Labor History, 50 (2009), pp.~347--350,
  \url{https://doi.org/10.1080/00236560903020930}.

\bibitem{hvarfner2024}
{\sc C.~Hvarfner, E.~O. Hellsten, and L.~Nardi}, {\em Vanilla {{Bayesian
  Optimization Performs Great}} in {{High Dimensions}}}, in Proceedings of the
  41st {{International Conference}} on {{Machine Learning}}, vol.~235 of
  Proceedings of {{Machine Learning Research}}, PMLR, 2024, pp.~20793--20817,
  \url{https://doi.org/10.48550/arxiv.2402.02229}.

\bibitem{ipsen1998}
{\sc I.~C.~F. Ipsen and C.~D. Meyer}, {\em The {{Idea Behind Krylov Methods}}},
  The American Mathematical Monthly, 105 (1998), pp.~889--899,
  \url{https://doi.org/10.1080/00029890.1998.12004985}.

\bibitem{jayaram2013}
{\sc B.~Jayaram and F.~Klawonn}, {\em Can {{Fuzzy Clustering Avoid Local
  Minima}} and {{Undesired Partitions}}?}, in Computational {{Intelligence}} in
  {{Intelligent Data Analysis}}, C.~Moewes and A.~N{\"u}rnberger, eds.,
  vol.~445, Springer Berlin Heidelberg, Berlin, Heidelberg, 2013, pp.~31--44,
  \url{https://doi.org/10.1007/978-3-642-32378-2_3}.

\bibitem{jin2009}
{\sc C.~Jin and X.~C. Cai}, {\em A preconditioned recycling {{GMRES}} solver
  for stochastic helmholtz problems}, Communications in Computational Physics,
  6 (2009), pp.~342--353, \url{https://doi.org/10.4208/cicp.2009.v6.p342}.

\bibitem{jones1998}
{\sc D.~R. Jones, M.~Schonlau, and W.~J. Welch}, {\em Efficient {{Global
  Optimization}} of {{Expensive}} {{Black-Box Functions}}}, Journal of Global
  Optimization, 13 (1998), pp.~455--492,
  \url{https://doi.org/10.1023/A:1008306431147}.

\bibitem{kalczynski2024}
{\sc P.~Kalczynski and Z.~Drezner}, {\em Further {{Analysis}} of the {{Weber
  Problem}}}, Networks and Spatial Economics,  (2024),
  \url{https://doi.org/10.1007/s11067-024-09627-1}.

\bibitem{keese2004}
{\sc A.~Keese}, {\em Numerical {{Solution}} of {{Systems}} with {{Stochastic
  Uncertainties}}: {{A General Purpose Framework}} for {{Stochastic Finite
  Elements}}}, PhD thesis, Universit\"atsbibliothek Braunschweig, Apr. 2004,
  \url{https://doi.org/10.24355/DBBS.084-200511080100-436}.

\bibitem{khoo2019switchnet}
{\sc Y.~Khoo and L.~Ying}, {\em {{SwitchNet}}: A neural network model for
  forward and inverse scattering problems}, SIAM Journal on Scientific
  Computing, 41 (2019), pp.~A3182--A3201.

\bibitem{lara2018}
{\sc C.~L. Lara, F.~Trespalacios, and I.~E. Grossmann}, {\em Global
  optimization algorithm for capacitated multi-facility continuous
  location-allocation problems}, Journal of Global Optimization, 71 (2018),
  pp.~871--889, \url{https://doi.org/10.1007/s10898-018-0621-6}.

\bibitem{lax1989}
{\sc P.~D. Lax and R.~S. Phillips}, {\em Scattering Theory}, no.~v. 26 in Pure
  and Applied Mathematics, Academic Press, Boston, rev. ed~ed., 1989.

\bibitem{lerer2024multigrid}
{\sc B.~Lerer, I.~{Ben-Yair}, and E.~Treister}, {\em Multigrid-augmented deep
  learning preconditioners for the {{Helmholtz}} equation using compact
  implicit layers}, SIAM Journal on Scientific Computing, 46 (2024),
  pp.~S123--S144.

\bibitem{luong2021}
{\sc P.~Luong, D.~Nguyen, S.~Gupta, S.~Rana, and S.~Venkatesh}, {\em Adaptive
  cost-aware {{Bayesian}} optimization}, Knowledge-Based Systems, 232 (2021),
  p.~107481, \url{https://doi.org/10.1016/j.knosys.2021.107481}.

\bibitem{luz2020}
{\sc I.~Luz, M.~Galun, H.~Maron, R.~Basri, and I.~Yavneh}, {\em Learning
  algebraic multigrid using graph neural networks}, in Proceedings of the 37th
  International Conference on Machine Learning, {{ICML}}'20, JMLR.org, 2020.

\bibitem{neal1997}
{\sc R.~M. Neal}, {\em Monte carlo implementation of gaussian process models
  for bayesian regression and classification}, 1997,
  \url{https://arxiv.org/abs/physics/9701026},
  \url{https://arxiv.org/abs/physics/9701026}.

\bibitem{nedelec2001}
{\sc J.-C. N{\'e}d{\'e}lec}, {\em Acoustic and {{Electromagnetic Equations}}},
  vol.~144 of Applied {{Mathematical Sciences}}, Springer New York, New York,
  NY, 2001, \url{https://doi.org/10.1007/978-1-4757-4393-7}.

\bibitem{parks2006}
{\sc M.~L. Parks, E.~{de Sturler}, G.~Mackey, D.~D. Johnson, and S.~Maiti},
  {\em Recycling {{Krylov Subspaces}} for {{Sequences}} of {{Linear Systems}}},
  SIAM Journal on Scientific Computing, 28 (2006), pp.~1651--1674,
  \url{https://doi.org/10.1137/040607277}.

\bibitem{pearson2020}
{\sc J.~W. Pearson and J.~Pestana}, {\em Preconditioners for {{Krylov}}
  subspace methods: {{An}} overview}, GAMM Mitteilungen, 43 (2020), pp.~1--35,
  \url{https://doi.org/10.1002/gamm.202000015}.

\bibitem{pellissetti2000}
{\sc M.~F. Pellissetti and R.~Ghanem}, {\em Iterative solution of systems of
  linear equations arising in the context of stochastic finite elements},
  Advances in Engineering Software, 31 (2000), pp.~607--616,
  \url{https://doi.org/10.1016/S0965-9978(00)00034-X}.

\bibitem{pembery2020}
{\sc O.~R. Pembery}, {\em The {{Helmholtz Equation}} in {{Heterogeneous}} and
  {{Random Media}}: {{Analysis}} and {{Numerics}}}, PhD thesis, University of
  Bath, 2020.

\bibitem{poulson2013}
{\sc J.~Poulson, B.~Engquist, S.~Li, and L.~Ying}, {\em A {{Parallel Sweeping
  Preconditioner}} for {{Heterogeneous 3D Helmholtz Equations}}}, SIAM Journal
  on Scientific Computing, 35 (2013), pp.~C194--C212,
  \url{https://doi.org/10.1137/120871985}.

\bibitem{powell2012}
{\sc C.~E. Powell and D.~J. Silvester}, {\em Preconditioning {{Steady-State
  Navier--Stokes Equations}} with {{Random Data}}}, SIAM Journal on Scientific
  Computing, 34 (2012), pp.~A2482--A2506,
  \url{https://doi.org/10.1137/120870578}.

\bibitem{rasmussen2000}
{\sc C.~E. Rasmussen and C.~K.~I. Williams}, {\em Gaussian {{Processes}} for
  {{Machine Learning}}}, vol.~7, The MIT Press, 2000.

\bibitem{rousseeuw1987}
{\sc P.~J. Rousseeuw}, {\em Silhouettes: {{A}} graphical aid to the
  interpretation and validation of cluster analysis}, Journal of Computational
  and Applied Mathematics, 20 (1987), pp.~53--65,
  \url{https://doi.org/10.1016/0377-0427(87)90125-7}.

\bibitem{ruspini2019}
{\sc E.~H. Ruspini, J.~C. Bezdek, and J.~M. Keller}, {\em Fuzzy {{Clustering}}:
  {{A Historical Perspective}}}, IEEE Computational Intelligence Magazine, 14
  (2019), pp.~45--55, \url{https://doi.org/10.1109/MCI.2018.2881643}.

\bibitem{scroggs2022}
{\sc M.~W. Scroggs, I.~A. Baratta, C.~N. Richardson, and G.~N. Wells}, {\em
  Basix: A runtime finite element basis evaluation library}, Journal of Open
  Source Software, 7 (2022), pp.~3982--3982,
  \url{https://doi.org/10.21105/joss.03982}.

\bibitem{scroggs2022a}
{\sc M.~W. Scroggs, J.~S. Dokken, C.~N. Richardson, and G.~N. Wells}, {\em
  Construction of {{Arbitrary Order Finite Element Degree-of-Freedom Maps}} on
  {{Polygonal}} and {{Polyhedral Cell Meshes}}}, ACM Transactions on
  Mathematical Software, 48 (2022), pp.~1--23,
  \url{https://doi.org/10.1145/3524456}.

\bibitem{shende2022}
{\sc S.~Shende, A.~Gillman, P.~Buskohl, and K.~Vemaganti}, {\em Systematic cost
  analysis of gradient- and anisotropy-enhanced {{Bayesian}} design
  optimization}, Structural and Multidisciplinary Optimization, 65 (2022),
  p.~235, \url{https://doi.org/10.1007/s00158-022-03324-8}.

\bibitem{sherali1988}
{\sc H.~D. Sherali and F.~L. Nordai}, {\em {{NP-Hard}}, {{Capacitated}},
  {{Balanced}} {\emph{p}} -{{Median Problems}} on a {{Chain Graph}} with a
  {{Continuum}} of {{Link Demands}}}, Mathematics of Operations Research, 13
  (1988), pp.~32--49, \url{https://doi.org/10.1287/moor.13.1.32}.

\bibitem{shirron1998}
{\sc J.~J. Shirron and I.~Babu{\v s}ka}, {\em A comparison of approximate
  boundary conditions and infinite element methods for exterior {{Helmholtz}}
  problems}, Computer Methods in Applied Mechanics and Engineering, 164 (1998),
  pp.~121--139, \url{https://doi.org/10.1016/S0045-7825(98)00050-4}.

\bibitem{snoek2012}
{\sc J.~Snoek, H.~Larochelle, and R.~P. Adams}, {\em Practical {{Bayesian
  Optimization}} of {{Machine Learning Algorithms}}}, in Advances in {{Neural
  Information Processing Systems}}, vol.~25, Curran Associates, Inc., 2012.

\bibitem{tartakovsky2006}
{\sc D.~M. Tartakovsky and D.~Xiu}, {\em Stochastic analysis of transport in
  tubes with rough walls}, Journal of Computational Physics, 217 (2006),
  pp.~248--259, \url{https://doi.org/10.1016/j.jcp.2006.02.029}.

\bibitem{thorndike1953}
{\sc R.~L. Thorndike}, {\em Who belongs in the family?}, Psychometrika, 18
  (1953), pp.~267--276, \url{https://doi.org/10.1007/BF02289263}.

\bibitem{tibshirani2001}
{\sc R.~Tibshirani, G.~Walther, and T.~Hastie}, {\em Estimating the {{Number}}
  of {{Clusters}} in a {{Data Set Via}} the {{Gap Statistic}}}, Journal of the
  Royal Statistical Society Series B: Statistical Methodology, 63 (2001),
  pp.~411--423, \url{https://doi.org/10.1111/1467-9868.00293}.

\bibitem{venkovic2023}
{\sc N.~Venkovic}, {\em Preconditioning Strategies for Stochastic Elliptic
  Partial Differential Equations}, PhD thesis, 2023.

\bibitem{wang2019}
{\sc G.~Wang and Q.~Liao}, {\em Efficient {{Spectral Stochastic Finite Element
  Methods}} for {{Helmholtz Equations}} with {{Random Inputs}}}, East Asian
  Journal on Applied Mathematics, 9 (2019), pp.~601--621,
  \url{https://doi.org/10.4208/eajam.140119.160219}.

\bibitem{wathen2015}
{\sc A.~J. Wathen}, {\em Preconditioning}, Acta Numerica, 24 (2015),
  pp.~329--376, \url{https://doi.org/10.1017/S0962492915000021}.

\bibitem{weiszfeld1937}
{\sc E.~Weiszfeld}, {\em Sur le point pour lequel la somme des distances de n
  points donnes est minimum}, Tohoku Mathematical Journal, 43 (1937),
  pp.~355--386.

\bibitem{xie2024}
{\sc Q.~Xie, R.~Astudillo, P.~Frazier, Z.~Scully, and A.~Terenin}, {\em
  Cost-aware {{Bayesian}} optimization via the {{Pandora}}'s {{Box Gittins}}
  index}, in 38th {{Conference}} on {{Neural Information Processing Systems}},
  arXiv, 2024, \url{https://arxiv.org/abs/2406.20062}.

\bibitem{xiu2006}
{\sc D.~Xiu and D.~M. Tartakovsky}, {\em Numerical methods for differential
  equations in random domains}, SIAM Journal on Scientific Computing, 28
  (2006), pp.~1167--1185, \url{https://doi.org/10.1137/040613160}.

\end{thebibliography}
    \clearpage

\end{document}